\documentclass[12pt]{article}
\usepackage{latexsym,amsfonts,amsthm,amsmath,amscd,amssymb}
\usepackage[dvips]{graphicx}

\newtheorem{lemma}{Lemma}[section]
\newtheorem{thm}[lemma]{Theorem}
\newtheorem{rem}[lemma]{Remark}
\newtheorem{prop}[lemma]{Proposition}
\newtheorem{cor}[lemma]{Corollary}

\newtheorem{defn}[lemma]{Definition}

\newcommand{\dimo}[1]{\vspace{2pt}\noindent\textit{Proof of \ref{#1}}.\ }
\newcommand{\finedimo}{{\hfill\hbox{$\square$}\vspace{2pt}}}

\newcommand\matS{{\mathbb{S}}}

\newcommand\matC{{\mathbb{C}}}

\newcommand\matK{{\mathbb{K}}}

\newcommand\calI{{\mathcal{I}}}

\renewcommand{\hbar}{{\overline{h}}}

\newfont{\Got}{eufm10 scaled 1200}

\newcommand{\mycap} [1] {\caption{\footnotesize{#1}}}

\newcommand{\Adm}{\mathop{\rm Adm}\nolimits}

\begin{document}

\title{On colored Turaev-Viro invariants\\ for links in arbitrary $3$-manifolds}

\author{Ekaterina~{\textsc Pervova}\thanks{The first author was supported by the
Marie Curie International Incoming Fellowship MIIF1-CT-2006-038734 within the
6th European Community Framework Programme} \and
\addtocounter{footnote}{5} Carlo~{\textsc Petronio}\thanks{The
second author gratefully thanks the University of Melbourne, the
Universit\'e Paul Sabatier in Toulouse, the Universit\'e Louis
Pasteur in Strasbourg and the Columbia University in New York
City, where part of this research was carried out} }

\maketitle

\begin{abstract}
\noindent We consider certain invariants of links in 3-manifolds,
obtained by a specialization of the Turaev-Viro invariants of
3-manifolds, that we call \emph{colored Turaev-Viro invariants}.
Their construction is based on a presentation of a pair $(M,L)$,
where $M$ is a closed oriented 3-manifold and $L\subset M$ is an
oriented link, by a triangulation of $M$ such that each component
of $L$ is an edge. We analyze some basic properties of these
invariants, including the behavior under connected sums of pairs
away and along links. These properties allow us to provide
examples of links in $\matS^3$ having the same HOMFLY polynomial
and the same Kauffman polynomial but distinct Turaev-Viro
invariants, and similar examples for the Alexander polynomial. We
also investigate the relations between the Turaev-Viro invariants
of $(M,L)$ and those of $M\setminus L$, showing that they are
sometimes but not always determined by each other.

\smallskip

\noindent MSC (2000): 57M27

\end{abstract}

\section*{Introduction}

The ideas of quantum topology \cite{Turaev3} have been
repeatedly applied
over the time to construct invariants of links in $3$-manifolds.
These constructions have led, among other things, to a computation of
the Jones polynomial of links in $\matS^3$ via a
2-dimensional~\cite{KirResh,Turaev2} or a
3-dimensional~\cite{TuraevIHES} state-sum model. In particular, in
the latter model the invariants are calculated as certain
state-sums on triangulations of link exteriors, with state-sums
based on the so-called quantum $6j$-symbols associated with the
quantum group $U_q(sl_2(\mathbb{C}))$~\cite{q-symbols,KirResh}.
The model furthermore takes into account a certain geometric
picture on the boundary of the exterior.

A different approach was used in~\cite{Barrett_et_al}, based on the
notion of a \emph{triangulation of a pair} $(M,L)$, namely a
triangulation of $M$ such that $L$ consists of edges, where $M$ is a
$3$-manifold and $L$ is a link in $M$. The invariants are defined as
state-sums on such triangulations, using the quantum $6j$-symbols
associated with $U_q(sl_2(\mathbb{C}))$ and restricting only to
states such that the edges contained in $L$ have a certain
pre-determined color. The main achievement of~\cite{Barrett_et_al} with
respect to these invariants was to establish a formula relating
their Fourier transforms to the so-called relativistic spin network
invariant, and to exhibit a relation with the
Witten-Reshetikhin-Turaev invariant of a certain framed link in a
$3$-manifold (where the underlying link is different from the initial one).

The aim of the present paper is to study in greater detail the
construction of~\cite{Barrett_et_al}, settling some of the issues
not addressed there. To this end we re-define the invariants (that
we call \emph{colored Turaev-Viro invariants}) on the basis of
algebraic initial data satisfying certain axioms (closely
following the style of~\cite{TV}). We note that our definition of
initial data differs from that in~\cite{TV} because of some
further requirements that we introduce. We then investigate the
behaviour of the invariants under the operations of connected sum
of pairs away and along links
(Theorems~\ref{ordinary:connected:thm} and~\ref{2-connected:thm}).
This point was not touched in~\cite{Barrett_et_al} even for the
specific invariants considered there.

We then proceed to consider two instances of initial data,
introduced in~\cite{TV} and
in~\cite{epsilon,Ma} respectively,
that satisfy our additional requirements. The aim here is to test the
effectiveness of the invariants by comparing them with other known
ones. In particular, our computations show that for links in
$\matS^3$ the colored Turaev-Viro invariants are independent of many
known ones, such as the HOMFLY polynomial, the Kauffman polynomial, and
the Alexander polynomial. Indeed, for each of these polynomials we
provide an example of pairs of links that are not distinguished by
the polynomial but are distinguished by one of the colored
Turaev-Viro invariants. Finally, we discuss relations with the
invariants of the link complements, showing that they are sometimes
but not always determined by each other.

On technical grounds, we note that we restrict to triangulations
such that each link component consists of one edge only.
Furthermore, for the whole paper we employ the dual viewpoint of
the so-called special spines~\cite{Ma}. More precisely, we present pairs
$(M,L)$, where $M$ is a closed oriented $3$-manifold and $L\subset M$
is an oriented link, by certain \emph{o-spines}, which are
oriented special skeleta of $M$ with some marked regions.

\section{Links, o-spines, and Turaev-Viro\\ invariants of $3$-manifolds}

In this section we recall some known definitions and facts that
will be used in the rest of the paper.

\paragraph{Links in $3$-manifolds} In this paper we consider pairs
$(M,L)$, where $M$ is an oriented closed $3$-manifold and $L$ is an
oriented link in $M$. For such pairs we consider the notions of
\emph{ordinary connected sum} and of \emph{link-connected sum}.

We say that $(M,L)$ is the ordinary connected sum of $(M_1,L_1)$ and
$(M_2,L_2)$ if it is homeomorphic to the pair obtained by removing
from each $M_i$ a small open ball $B_i$ disjoint from $L_i$ and by
identifying $\partial B_1$ with $\partial B_2$ via an
orientation-reversing homeomorphism. We say that $(M,L)$ is the
link-connected sum of $(M_1,L_1)$ and $(M_2,L_2)$ if it is
homeomorphic to the pair obtained by removing from $M_i$ a small
open ball $B_i$ such that $B_i\cap L_i$ consists of one simple
unknotted arc and identifying $\partial B_1$ with $\partial B_2$ via
an orientation-reversing homeomorphism such that $\partial B_1\cap
L_1$ is identified with $\partial B_2\cap L_2$ in the only way that
induces a consistent orientation on $L$.

\begin{rem}
\emph{For connected $M_1$ and $M_2$ the ordinary connected sum is
well-defined and the link-connected sum depends only on the link
components involved, so it is well-defined for knots.}
\end{rem}

\paragraph{O-spines} The notion of o-spine extends the notion of
a special spine of a $3$-manifold introduced in~\cite{Casler}.

\begin{figure}[h]
    \begin{center}
    \includegraphics[scale=0.6]{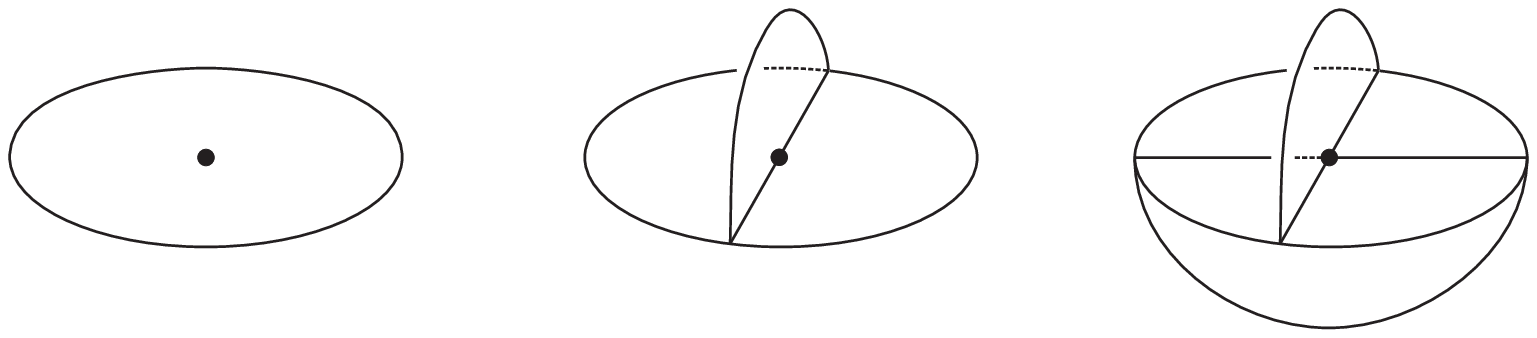}
    \mycap{Type I, II, and III neighborhoods of points in a special polyhedron}
    \label{neighborhoods:fig}
    \end{center}
    \end{figure}

\begin{defn}
\emph{A polyhedron $P$ is called \emph{special} if}
\begin{itemize}
\item \emph{each point of $P$ has neighbourhood of type I, II, or
III as shown in Fig.~\ref{neighborhoods:fig};}
\item \emph{the connected components of the set of type I points
are open discs (henceforth called \emph{regions});}
\item \emph{the connected components of the set of type II points are
open intervals (henceforth called \emph{edges}).}
\end{itemize}
\end{defn}

Type III points are called \emph{vertices}. The union of the edges and
the vertices is called the \emph{singular graph} of $P$ and is denoted
by $S(P)$. We also denote the set of all vertices by $V(P)$, the set
of all edges by $E(P)$, and the set of all regions by $R(P)$.

A subpolyhedron of a closed $3$-manifold $M$ is called a
\emph{spine} of $M$ if its complement is homeomorphic to an open
3-ball. Two manifolds with homeomorphic spines need not be
homeomorphic but they are if the spines are special~\cite{Casler}.
For this and more results on special spines see~\cite{Ma}.

\begin{figure}
    \begin{center}
    \includegraphics[scale=0.6]{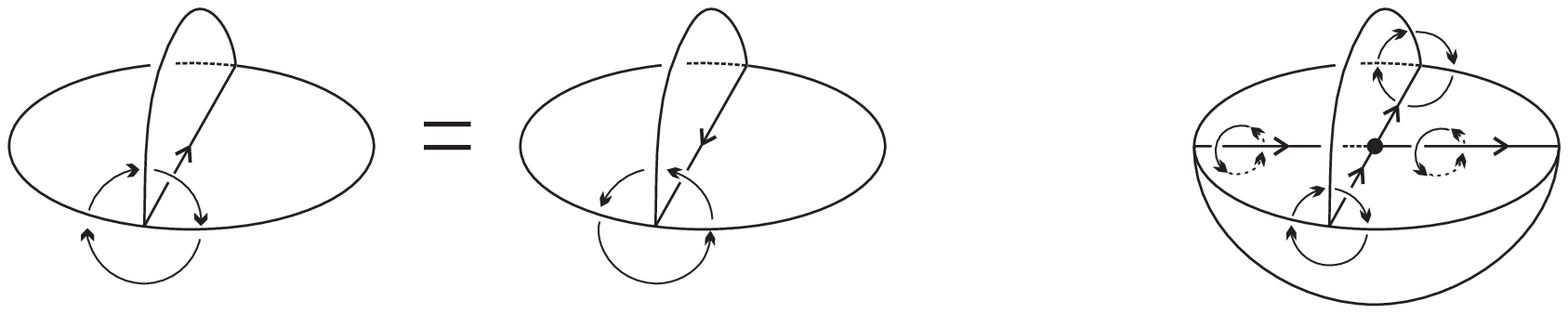}
    \mycap{Convention on screw-orientations and compatibility at vertices
    for an oriented special polyhedron}
    \label{orientation-spine:fig}
    \end{center}
    \end{figure}

An \emph{orientation} on a special polyhedron, as defined
in~\cite{BePe}, is a choice of a ``screw-orientation'' along the
edges (as in Fig.~\ref{orientation-spine:fig}-left), with
compatibility at vertices as in
Fig.~\ref{orientation-spine:fig}-right.
If $P$ is a special spine of an orientable $M$, an orientation on
$M$ determines one on $P$, and, conversely, an orientation on $P$
uniquely determines an orientation on $M$.

We define a \emph{marking} on a special polyhedron $P$ to be the
selection of some regions of $P$ together with the choice of an
orientation for each selected region.

\begin{defn}\label{o-spine:def}
\emph{If $M$ is a closed oriented $3$-manifold and $L\subset M$ is an
oriented link then a marked oriented special subpolyhedron $P$ of
$M$ is called an \emph{o-spine} of $(M,L)$ if:
\begin{enumerate}
 \item $P$ and $M$ are consistently oriented;
 \item $L$ intersects $P$ away from $S(P)$ and transversely;
 \item\label{o-spine:compl:cond} Each connected component $B$ of $M\setminus P$ is a
3-ball and $B\cap L$ is either empty or an unknotted proper arc in $B$;
 \item Each component of $L$ intersects $P$ at one point only;
 \item The marked regions of $P$ are precisely those intersected by $L$,
and they are oriented consistently with the orientations of $L$ and
$M$.
\end{enumerate}}
\end{defn}
Since each component of $L$ intersects only one (marked) region of
$P$, we will call this region \emph{dual} to the component.

\begin{figure}
    \begin{center}
    \includegraphics[scale=0.6]{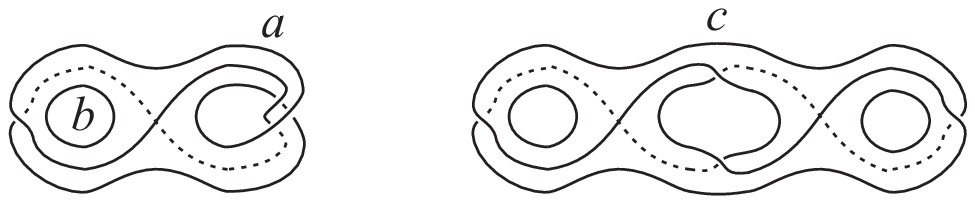}
    \mycap{The abalone and Bing's House }
    \label{Abalone-Bing:fig}
    \end{center}
    \end{figure}

We will always represent a special polyhedron $P$ by a picture
suggesting how a regular neighbourhood of $S(P)$ is embedded in
$\matS^3$, and we will endow $P$ with the orientation induced by
this embedding. Two examples are shown in
Fig.~\ref{Abalone-Bing:fig}.
In both cases it is not hard to check that $P$ is a special spine of $\matS^3$. The spine
in Fig.~\ref{Abalone-Bing:fig}-left is called the abalone and it has two
regions $a$ and $b$. If we select the region $a$ and we give it any
orientation, the resulting marked spine is an o-spine of the unknot,
whereas if we select $b$ we get the right-handed trefoil (again the
orientation of $b$ is immaterial).
The spine in Fig.~\ref{Abalone-Bing:fig}-right is called the Bing house.
It has three regions, marking any of which we get an o-spine of the
unknot. The following fact is known \cite{Amendola,Pe1}.

\begin{prop}
Every pair $(M,L)$ as above admits an o-spine.
\end{prop}

\begin{figure}
    \begin{center}
    \includegraphics[scale=0.7]{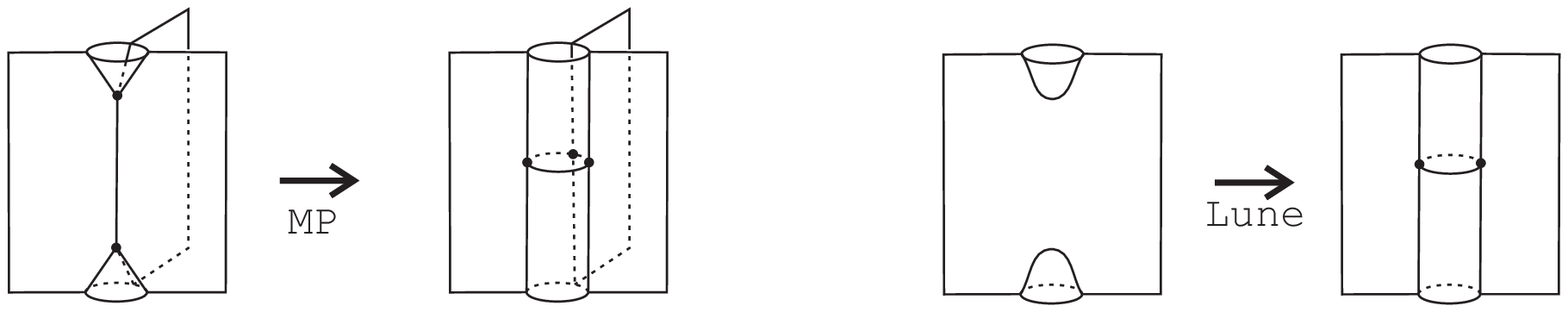}
    \mycap{Moves on oriented special polyhedra, with local orientations
    induced by the embedding in $\matS^3$}
    \label{moves:fig}
    \end{center}
    \end{figure}

\paragraph{Moves on o-spines}
Consider the moves $\texttt{MP}$ and $\texttt{Lune}$ shown in
Fig.~\ref{moves:fig}.
Of course if we apply $\texttt{MP}$ (which stands for
``Matveev-Piergallini,'' see~\cite{Ma} for the appropriate
references), or its inverse, to a special spine of a $3$-manifold
we get a new special spine of the same manifold. Moreover, it is
known~\cite{Ma} that any two special spines of the same
$3$-manifold can be transformed into each other by a sequence of
$\texttt{MP}$ moves and their inverses provided that they have at
least two vertices. Adding the \texttt{Lune} move the same fact
holds without restrictions on the number of vertices.

The positive move $\texttt{MP}$ has an obvious extension to the
context of \emph{marked} spines, because all old regions survive
and we can define the new region \emph{not} to be marked. This
applies to the inverse of $\texttt{MP}$ only if the triangular
region that disappears with the move is \emph{not} marked. Such an
inverse $\texttt{MP}$-move is called \emph{admissible}. It is then
known \cite{Amendola} that any two o-spines of the same pair can
be transformed into each other by a sequence of admissible
$\texttt{MP}$-moves and their inverses provided that they have at
least two vertices each.

The  $\texttt{Lune}$ move (see Fig.~\ref{moves:fig}-right) has
similar properties: whenever applied to a special spine of a
manifold or to an o-spine of a pair, it yields a special spine of
the same manifold or, respectively, an o-spine of the same pair.
Application of an inverse $\texttt{Lune}$-move may lead to a
non-special polyhedron and, just as an inverse $\texttt{MP}$-move,
to the disappearance of a marked region. An inverse
$\texttt{Lune}$-move is called \emph{admissible} if neither of
these happens. Any two o-spines of the same pair $(M,L)$,
regardless of the number of vertices, are related by a sequence of
admissible \texttt{MP} and \texttt{Lune} moves and their inverses.

\paragraph{Turaev-Viro initial data} Let $\matK$
be a commutative ring with
unity. We recall that an initial datum for a Turaev-Viro
invariant~\cite{TV} consists of a finite set $I$, a function that
assigns to each $i\in I$ its \emph{weight} $w_i\in \matK^*$, a certain
distinguished element $w\in \matK^*$, a distinguished set of unordered
triples of elements of $I$, that are called \emph{admissible}, and a
function that assigns to each so-called \emph{admissible $6$-tuple}
its \emph{$6j$-symbol}, which is an element of $\matK$.

More precisely, an ordered 6-tuple $(i,j,k,l,m,n)$ of elements of
$I$ is \emph{admissible} if the unordered triples $(i,j,k)$,
$(i,m,n)$, $(j,l,n)$, $(k,l,m)$ are admissible. The $6j$-symbol
associated to this 6-tuple is denoted by
$$\left|\begin{array}{ccc}i & j & k \\ l & m & n\end{array}\right|.$$
The symbols must have the following symmetries:
$$\left.\begin{array}{ccccccl}
 & \left|\begin{array}{ccc}i & j & k \\ l & m & n\end{array}\right| &
 = &
\left|\begin{array}{ccc}j & i & k \\ m & l & n\end{array}\right| & =
& \left|\begin{array}{ccc}i & k & j \\ l & n & m\end{array}\right| & \\
 & & & & & & \\
 = & \left|\begin{array}{ccc}i & m & n \\ l & j & k\end{array}\right|
& = & \left|\begin{array}{ccc}l & m & k \\ i & j &
n\end{array}\right| & = & \left|\begin{array}{ccc}l & j & n \\ i & m
& k\end{array}\right| & .\end{array}\right.
$$

Let us define some properties that an initial datum may (or not)
satisfy. Notation and terminology come from the literature (see, for
instance,~\cite{Turaev3} and references therein).
\begin{itemize}
\item[$(*)$] \emph{Orthogonality}: For any $i',i'',j,k,l,n\in I$ such that the
triples $(i',j,k)$, $(i',l,n)$, $(i'',j,k)$, $(i'',l,n)$ are
admissible we have:
$$\sum\limits_mw_m\left|\begin{array}{ccc}i' & j & k \\ m & l & n\end{array}\right|
\left|\begin{array}{ccc}i'' & j & k \\ m & l & n\end{array}\right|=
w_i^{-1}\delta_{i'i''}$$ with sum taken over all $m$ such that both
symbols are defined;

\item[$(**)$]
\emph{Biedenharn-Elliot identity}: For any admissible
6-tuples\\ $(k,k',k'',j,j',j'')$ and $(k,k',k'',i,i',i'')$ we have:
$$\left.\begin{array}{rl}
 & \left|\begin{array}{ccc}k & k' & k'' \\ j & j' & j''\end{array}\right|\cdot
 \left|\begin{array}{ccc}k & k' & k'' \\ i & i' & i''\end{array}\right|\\
 & \\
 = & \sum\limits_nw_n\left|\begin{array}{ccc} j' & k & j'' \\ i'' & n & i'\end{array}\right|\cdot
 \left|\begin{array}{ccc}j & k' & j'' \\ i'' & n & i\end{array}\right|\cdot
 \left|\begin{array}{ccc} j & k'' & j' \\ i' & n & i\end{array}\right|
 \end{array}\right.
$$
with sum taken over all $n$'s such that the involved symbols are
defined;

\item[$(*\!*\!*)$] For any
$i\in I$ we have
$$\sum\limits_{j,k}w_jw_k=w\cdot w_i^{-1} $$
with sum taken over all $j,k$ such that $(i,j,k)$ is
admissible;

\item \emph{Irreducibility}: For any
$j,k\in I$ there exists a sequence $l_1,l_2,\ldots,l_n\in I$ with
$j=l_1$ and $k=l_n$ such that $(l_i,l_{i+1},l_{i+2})$ is admissible.
\end{itemize}

\paragraph{Turaev-Viro invariants} To any irreducible initial datum
satisfying Conditions~$(*)$, $(*\,*)$, and $(*\!*\!*)$ there corresponds a
$3$-manifold invariant of Turaev-Viro type~\cite[Section~1.3]{TV},
whose value on a given $3$-manifold $M$ is obtained as follows. If $P$
is a special spine of $M$, we call a \emph{coloring} of $P$ any map
$\eta:R(P)\rightarrow I$. A coloring is called
\emph{admissible} if for any edge $e$ of $P$ the colors of the three
germs of regions incident to $e$ form an admissible triple.

\begin{figure}
    \begin{center}
    \includegraphics[scale=0.6]{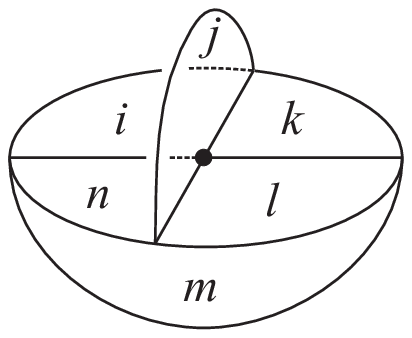}
    \mycap{A coloring around a vertex yields a $6j$-symbol}
    \label{vertex-symbol:fig}
    \end{center}
    \end{figure}

Denote the set of all admissible colorings of $P$ by $\Adm(P)$. If
$\eta\in\Adm(P)$ and $v\in V(P)$, we define $s_{\eta}(v)$ as the
$6j$-symbol
$$\left|\begin{array}{ccc}i & j & k \\ l & m &
n\end{array}\right|$$ where $i,j,k,l,m,n$ are the colors of the
germs of regions of $P$ incident to $v$, as shown in
Fig.~\ref{vertex-symbol:fig}. Note that ${s}_{\eta}(v)$ is
well-defined because of the symmetries of the $6j$-symbols.

We define the \emph{weight} of $\eta\in\Adm(P)$ as
$$\omega(\eta)=\prod_{v\in V(P)}{s}_{\eta}(v)\prod_{c\in R(P)}w_{\eta(c)}$$ and the
\emph{invariant of Turaev-Viro type} corresponding to the initial
datum as
$$TV(M)=w^{-1}\cdot\sum\limits_{\eta\in\Adm(P)}\omega(\eta).$$
As shown in~\cite{TV}, independence of the right-hand side of this
equality from the particular choice of $P$ follows from
Conditions~$(*)$ and~$(*\,*)$, which ensure invariance under $\texttt{Lune}$- and $\texttt{MP}$-moves
and their inverses. Condition~$(***)$ allows one to compute
the invariant using a slightly more general notion of special spine.
Namely, one requires $M\setminus P$ to be the union of some number
$m\geqslant 1$ of balls, not necessarily a single one, and one
replaces $w^{-1}$ in the formula defining $TV(M)$ by $w^{-m}$.

\begin{rem}\label{invt:for:bry:rem}
\emph{The definition of the invariant also makes sense in the context of manifolds
with boundary. Specifically, if $P$ is a special skeleton of a compact $3$-manifold
with boundary $M$ (\emph{i.e.} if $P$ is special and
$M\setminus P$ consists of the regular neighbourhood of the boundary and
of some number $m$ of open balls) then the invariant of Turaev-Viro type of $M$ is again
defined as $TV(M)=w^{-m}\sum\limits_{\eta\in\Adm(P)}\omega(\eta)$. Notice
also that this definition is different from that of the relative
Turaev-Viro invariant $\Omega_M(\alpha)$ introduced in~\cite[Section~1.4]{TV}.}
\end{rem}

\section{Turaev-Viro invariants of links}

In this section we describe certain additional requirements on the
initial data and state the definition of colored Turaev-Viro invariants of links.

\paragraph{Good initial data} The further conditions
on the initial datum that we will need are the following ones.
\begin{itemize}
\item \emph{Strong irreducibility}:
\begin{enumerate}
    \item[(i)] For any $i,j\in I$ there exists $k\in I$ such that the
    triple $(i,j,k)$ is admissible;
    \item[(ii)] For any two admissible triples of the form $(i,j,k)$ and $(i,m,n)$
    there exist $l,l'\in I$ such that the triples $(j,m,l)$,
    $(k,n,l)$, $(j,n,l')$, and $(k,m,l')$ are admissible;
    \item[(iii)] For any two admissible 6-tuples of the form  $(i,j,k,i_1,k_1,j_1)$
and $(i,j,k,i_2,k_2,j_2)$ there exists $l$ such that all the triples
    $(i_1,i_2,l)$, $(j_1,j_2,l)$, $(k_1,k_2,l)$ are admissible.
\end{enumerate}

\item [$(***\,*)$]
$\sum\limits_kw_k=w_iw_j \mbox{ for any }i,j\in I$, with sum taken over all $k$'s such that
the triple $(i,j,k)$ is admissible.
\end{itemize}

An initial datum $\mathcal{D}$ is called \emph{good} if it is
strongly irreducible and it satisfies Conditions~$(*)$, $(*\,*)$,
$(***)$, and~$(***\,*)$.

\paragraph{Turaev-Viro invariants of colored links} Let us fix a good
initial datum. We will now define an invariant $\calI_{\xi}(M,L)$ of
triples $(M,L,\xi)$, where $\xi$ is a coloring of the link components by
elements of $I$. Denote the components of $L$ by $L_1$, $\ldots$,
$L_n$, so $\xi(L_j)\in I$ is the color of $L_j$. Let $P$ be an
o-spine of $(M,L)$, let $\alpha_j$ be the region of $P$ dual to
$L_j$, and let $m$ be the number of complementary
balls having empty intersection with $L$. Denote by $\Adm_{\xi}(P)$
the set of all admissible colorings $\eta$ of $P$ such that
$\eta(\alpha_j)=\xi(L_j)$ for all $j$. Then
we set:
$$\calI_{\xi}(P)=w^{-m}\sum\limits_{\eta\in\Adm_{\xi}(P)}\omega(\eta).$$
The next result is proved basically by the same argument as in~\cite{TV}
(see also~\cite{Ma}). See Section~\ref{good:data:section} for details.

\begin{prop}\label{invt:well-defined:prop}
Let $P$ and $P'$ be o-spines of a pair $(M,L)$. Then
$\mathcal{I}_{\xi}(P)=\mathcal{I}_{\xi}(P')$.
\end{prop}

This proposition readily implies that we have a well-defined invariant
$\mathcal{I}_{\xi}(M,L)$ of triples $(M,L,\xi)$ computed as
$\calI_{\xi}(P)$ whenever $P$ is an o-spine of $(M,L)$. It also
follows from the definition that
$\calI(M,\emptyset)=TV(M)$ ---there is no coloring $\xi$
to choose on the empty link.

\begin{rem}
\emph{As it will follow from
the proof of Proposition~\ref{invt:well-defined:prop},
conditions~$(*)$,~$(*\,*)$, and~$(***)$ already imply that
$\calI_{\xi}(M,L)$ is well-defined. The other conditions imposed on
good initial data will be used to analyze how $\calI_{\xi}$ behaves
under connected sums.}
\end{rem}

\begin{rem}
\emph{For the Turaev-Viro initial data~\cite[Section~7]{TV} the
invariant $\calI_{\xi}$ was already defined
in~\cite{Barrett_et_al}.}
\end{rem}

\paragraph{The invariants $\calI_i$} Fixing a color $i\in I$ we can define
an invariant $\calI_i(M,L)$ by coloring with $i$ all the components
of $L$. When $M=\matS^3$, we denote the invariant by $\calI_i(L)$ only.

\section{General constructions}\label{gen:contr:sec}

In this section we describe several constructions that will be
used extensively throughout the paper.

\begin{figure}
    \begin{center}
    \includegraphics[scale=0.5]{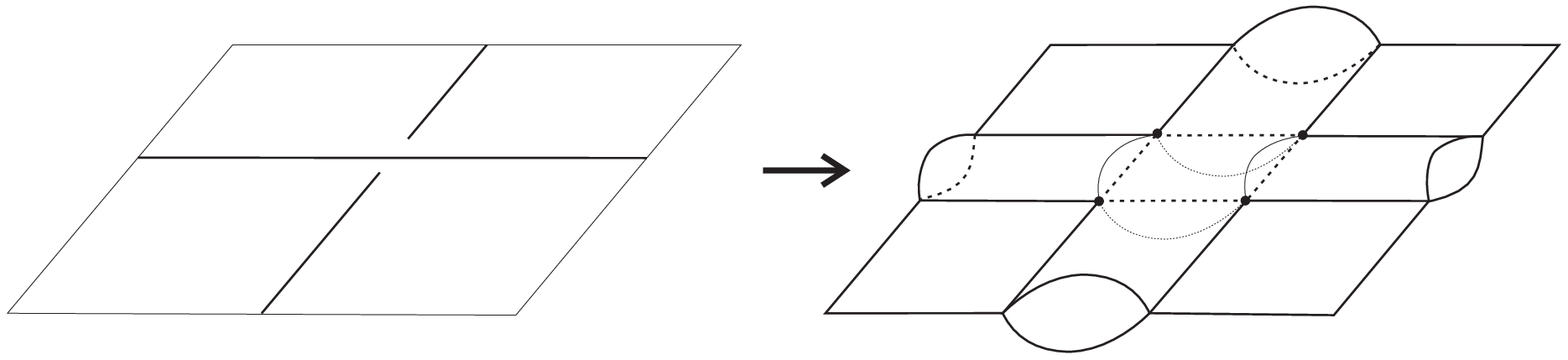}
    \mycap{Digging the tunnel in the neighbourhood of a crossing point}
    \label{crossing-point:fig}
    \end{center}
    \end{figure}

\begin{figure}
    \begin{center}
    \includegraphics[scale=0.5]{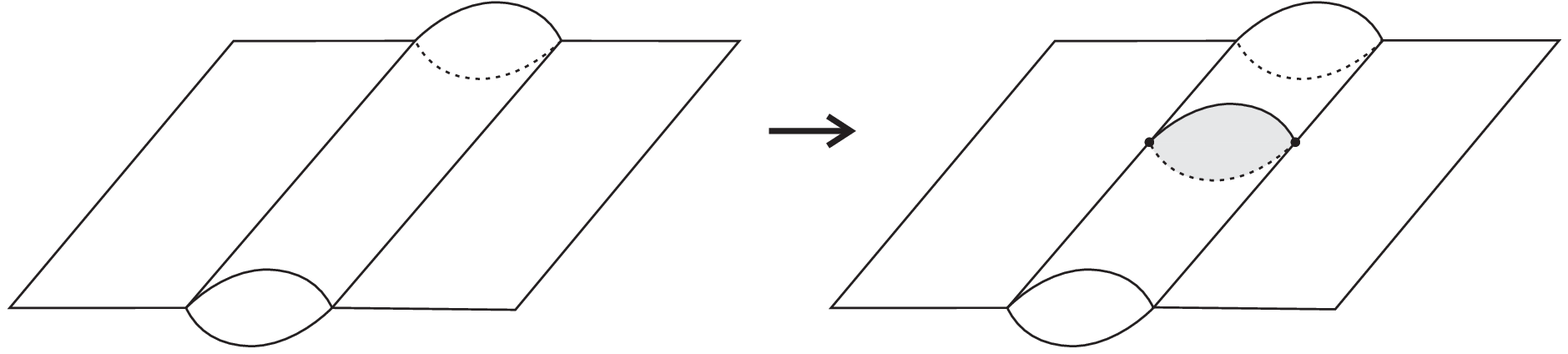}
    \mycap{Gluing a transverse disc}
    \label{transverse-disc:fig}
    \end{center}
    \end{figure}

\begin{figure}
    \begin{center}
    \includegraphics[scale=0.5]{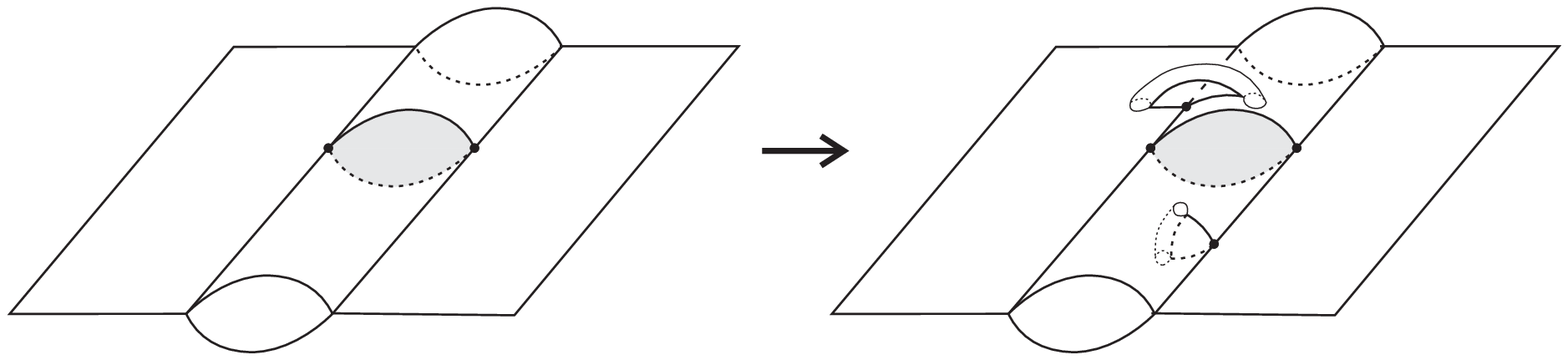}
    \mycap{Getting rid of the extra balls in the complement}
    \label{1-move-proj-spine:fig}
    \end{center}
    \end{figure}

\paragraph{An o-spine $P_{\Pi}(L)$ of $L\subset\matS^3$ from a
projection} Given a connected projection $\Pi$ of a link
$L\subset\matS^3$ on a 2-sphere $\matS^2\subset\matS^3$, we can
always construct an o-spine of $L$ as follows. For each component of
the link we ``dig a tunnel'' in $\matS^2$ along the projection.
Fig.~\ref{crossing-point:fig} shows the way to do this in the
neighbourhood of a crossing point.
We get an embedded special polyhedron $P'$ whose complement in
$\matS^3$ consists of two balls and a tubular neighbourhood $U(L)$ of the
link. For each connected component of $L$ we now add a
2-disc to $P'$ as in Fig.~\ref{transverse-disc:fig}, denoting the
resulting polyhedron by $P''$.
Since $\matS^3\setminus P'$ consists, besides $U(L)$, of two balls,
$P''$ is an o-spine of $(\matS^3,L)$: the marked discs are
those of Fig.~\ref{transverse-disc:fig}, and the orientation matters
are easily settled. We then perform on $P''$ the move of Fig.~\ref{1-move-proj-spine:fig}
(on a single component of $L$);
the result is still an o-spine of $(\matS^3,L)$, but now the complement of $P''$
is just a regular neighborhood of $L$, without additional balls.
We denote any polyhedron arising from this procedure by
$P_{\Pi}(L)$.

\begin{figure}
    \begin{center}
    \includegraphics[scale=0.5]{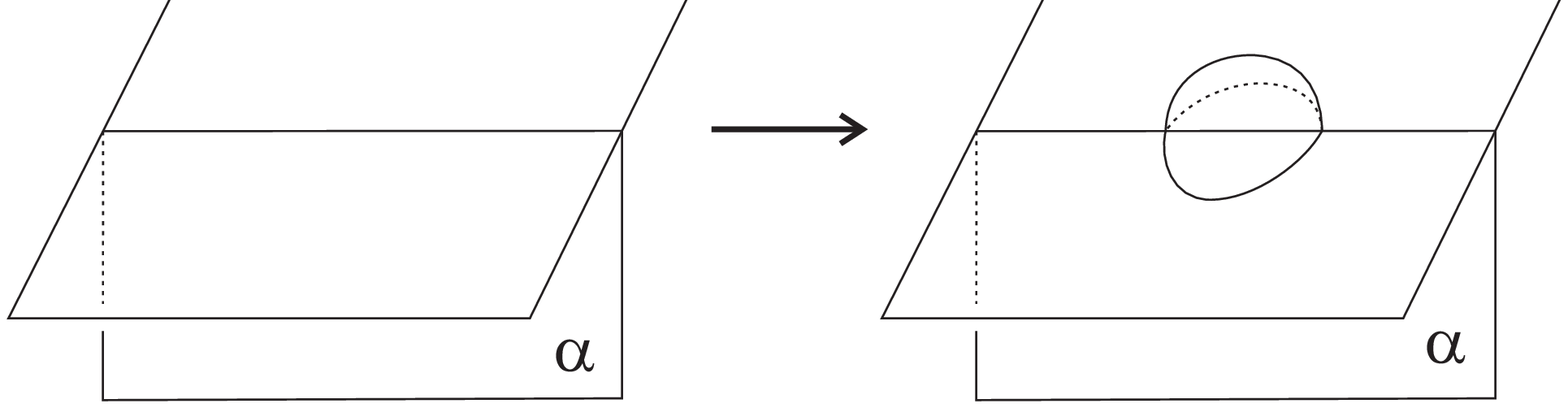}
    \mycap{Adding a bubble }
    \label{bubble:fig}
    \end{center}
    \end{figure}

\paragraph{A bubbled o-spine} Let $P$ be
a special polyhedron, and let $\alpha$ be a germ of region at
an edge of $P$. A move as in Fig.~\ref{bubble:fig} will be called
\emph{addition of a bubble}. The resulting polyhedron will be
denoted by $\mbox{{\tt Bub}}_{\alpha}(P)$.

\begin{figure}
    \begin{center}
    \includegraphics[scale=0.4]{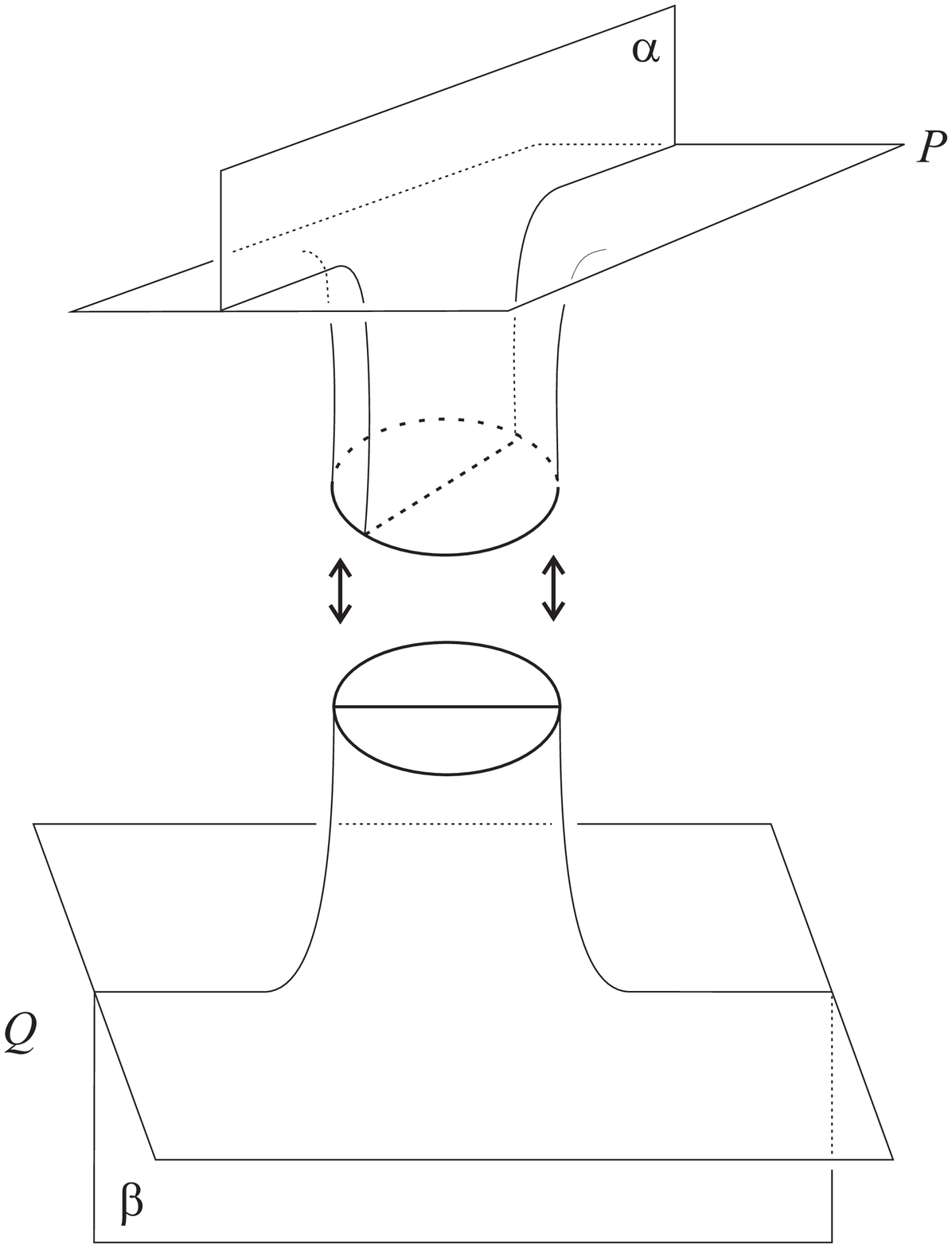}
    \mycap{Construction of $\mbox{{\tt
Join}}\big((P,\alpha),(Q,\beta)\big)$ }
    \label{join:fig}
    \end{center}
    \end{figure}

\paragraph{Join of polyhedra} Let $P$,
$Q$ be special polyhedra, let $\alpha$ be a germ of a region of $P$,
and let $\beta$ be a germ of a region of $Q$. We then denote by
$\mbox{{\tt Join}}\big((P,\alpha),(Q,\beta)\big)$ the special
polyhedron which is obtained as in Fig.~\ref{join:fig}.

\paragraph{O-spine of an
ordinary connected sum} Let $P$ be an o-spine of a pair $(M_1,L_1)$,
let $Q$ be an o-spine of $(M_2,L_2)$, and let $\alpha$ and $\beta$
be germs of regions of $P$ and of $Q$ respectively. Construct
$\mbox{{\tt Bub}}_{\alpha}(P)$ first. We remark that $\alpha$ gives
rise to three new germs in $\mbox{{\tt Bub}}_{\alpha}(P)$, exactly
one of which is incident to both the two new vertices (this is
the germ with boundary ``inside the bubble'', see
Fig.~\ref{bubble:fig}). We denote this germ by $\overline{\alpha}$ and
set
$$\mbox{{\tt OSum}}_{\alpha,\beta}(P,Q):=\mbox{{\tt
Join}}\big((\mbox{{\tt
Bub}}_{\alpha}(P),\overline{\alpha}),(Q,\beta)\big).$$

\begin{lemma}\label{construct:ordin:spine:lem}
The polyhedron $\mbox{{\tt OSum}}_{\alpha,\beta}(P,Q)$ is an o-spine
of $(M_1\#_oM_2,L_1\sqcup L_2)$.
\end{lemma}

\begin{figure}
    \begin{center}
    \includegraphics[scale=0.4]{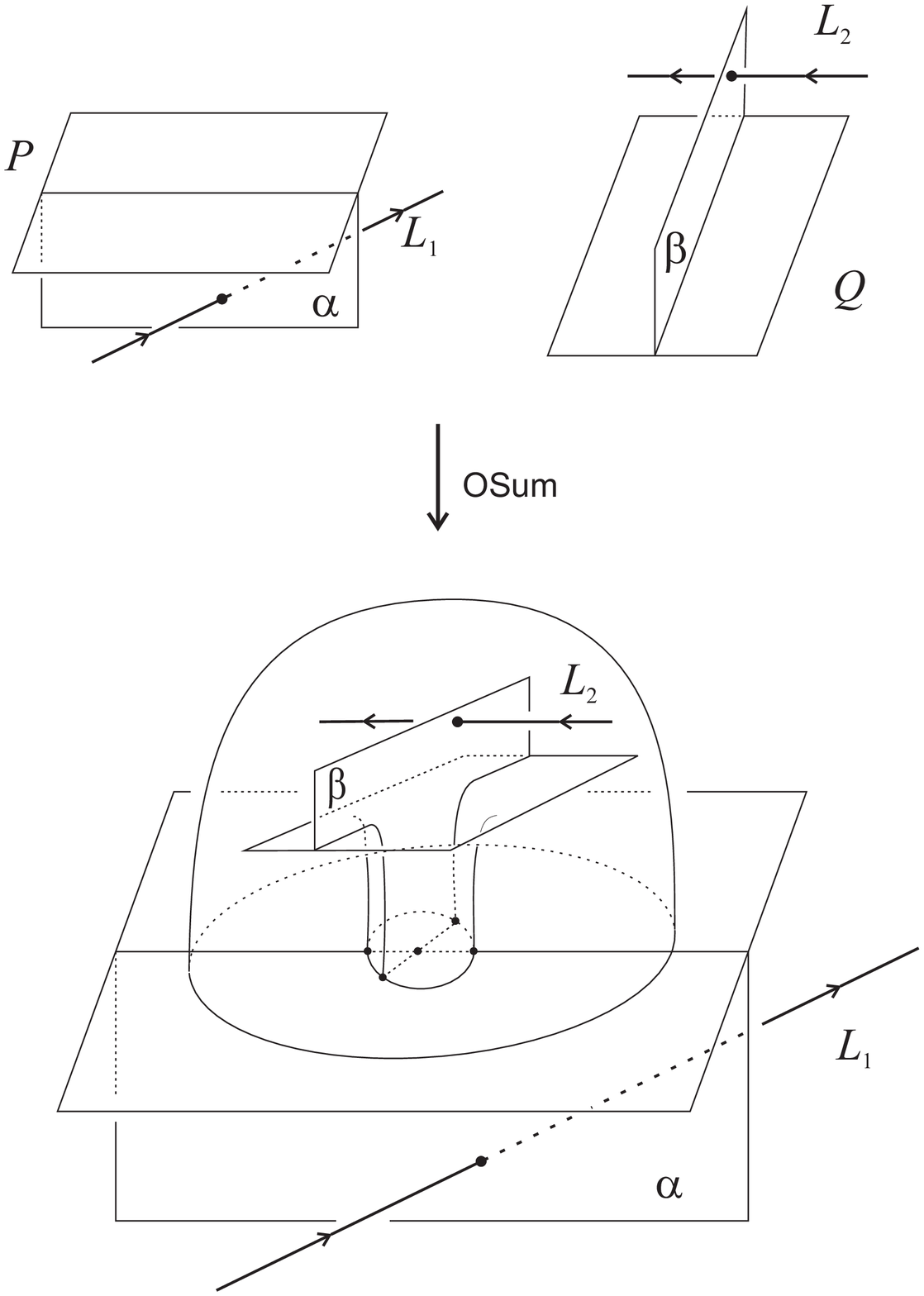}
    \mycap{Construction of $\mbox{{\tt OSum}}_{\alpha,\beta}(P,Q)$ }
    \label{ord-sum-spine:fig}
    \end{center}
    \end{figure}

The proof of this lemma is given by Fig.~\ref{ord-sum-spine:fig} (note that there may or not
be components of $L_1$ and $L_2$ dual to $\alpha$ and $\beta$).

\begin{figure}
    \begin{center}
    \includegraphics[scale=0.5]{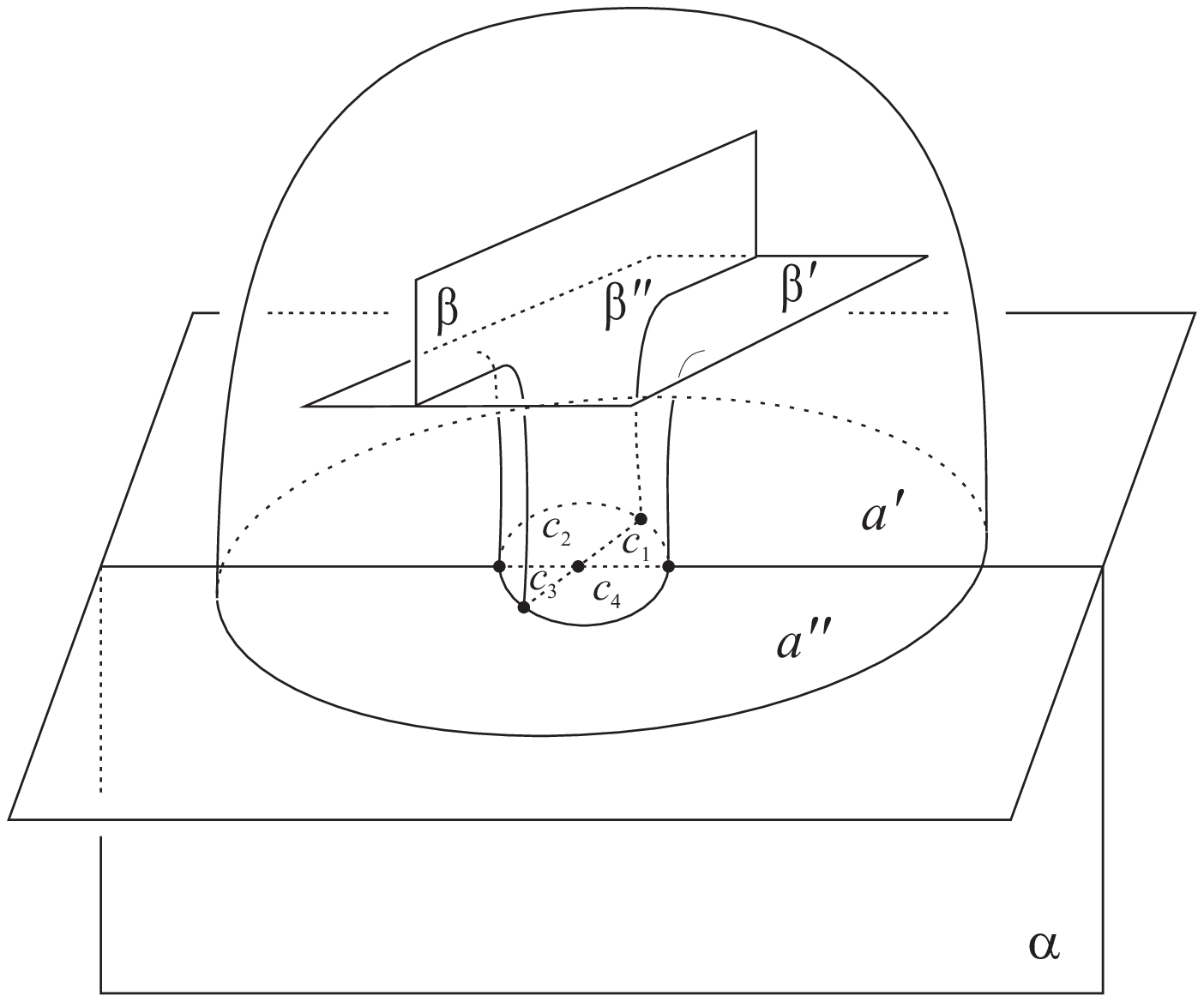}
    \mycap{
Notation for the regions of $\mbox{{\tt Join}}\big((\mbox{{\tt
Bub}}_{\alpha}(P),\overline{\alpha}),(Q,\beta)\big)$ used in the
construction of $\mbox{{\tt LSum}}_{\alpha,\beta}(P,Q)$ and below
in Section~\ref{good:data:section}}
    \label{join-cell-notation:fig}
    \end{center}
    \end{figure}

\paragraph{O-spine of a link-connected
sum} Let $P$ and $Q$ be o-spines of pairs $(M_1,L_1)$ and
$(M_2,L_2)$. Let $K_i$ be a component of $L_i$, with $K_1$ dual to
the region $\alpha$ of $P$ and $K_2$ dual to the region $\beta$ of
$Q$. We construct first the marked oriented polyhedron $\mbox{{\tt
Join}}\big((\mbox{{\tt
Bub}}_{\alpha}(P),\overline{\alpha}),(Q,\beta)\big)$ and denote by $c_1$,
$\ldots$, $c_4$ the four triangular regions that arise from
identifying two discs in $\mbox{{\tt Bub}}_{\alpha}(P)$ and $Q$ (see
Fig.~\ref{join-cell-notation:fig}).

For $i=1,2,3,4$ we denote by $\mu_\alpha(c_i)$ (respectively
$\mu_{\beta}(c_i)$) the orientation of $c_i$ such that $\alpha$
(respectively $\beta$) and $c_i$ induce opposite orientations on
their common edge. Then there are exactly two indices
$s_1,s_2\in\{1,2,3,4\}$ such that
$\mu_{\alpha}(c_{s_j})=\mu_{\beta}(c_{s_j})$. We then define the
polyhedron $\mbox{{\tt LSum}}_{\alpha,\beta}(P,Q)$ as the o-spine
obtained from $\mbox{{\tt Join}}\big((\mbox{{\tt
Bub}}_{\alpha}(P),\overline{\alpha}),(Q,\beta)\big)$ by removing
the markings from $\alpha$ and $\beta$, placing the marking
$\mu_{\alpha}(c_{i_1})=\mu_{\beta}(c_{s_1})$ on $c_{s_1}$, and
removing the region
$c_{s_2}$.

\begin{lemma}\label{construct:2-con:spine:lem}
Let $(M,L)$ be obtained by the link-connected sum of $(M_1,L_1)$ and
$(M_2,L_2)$ along $K_1$ and $K_2$. Then any $\mbox{{\tt
LSum}}_{\alpha,\beta}(P,Q)$ constructed as above is an o-spine of
$(M,L)$.
\end{lemma}

\begin{figure}
    \begin{center}
    \includegraphics[scale=0.4]{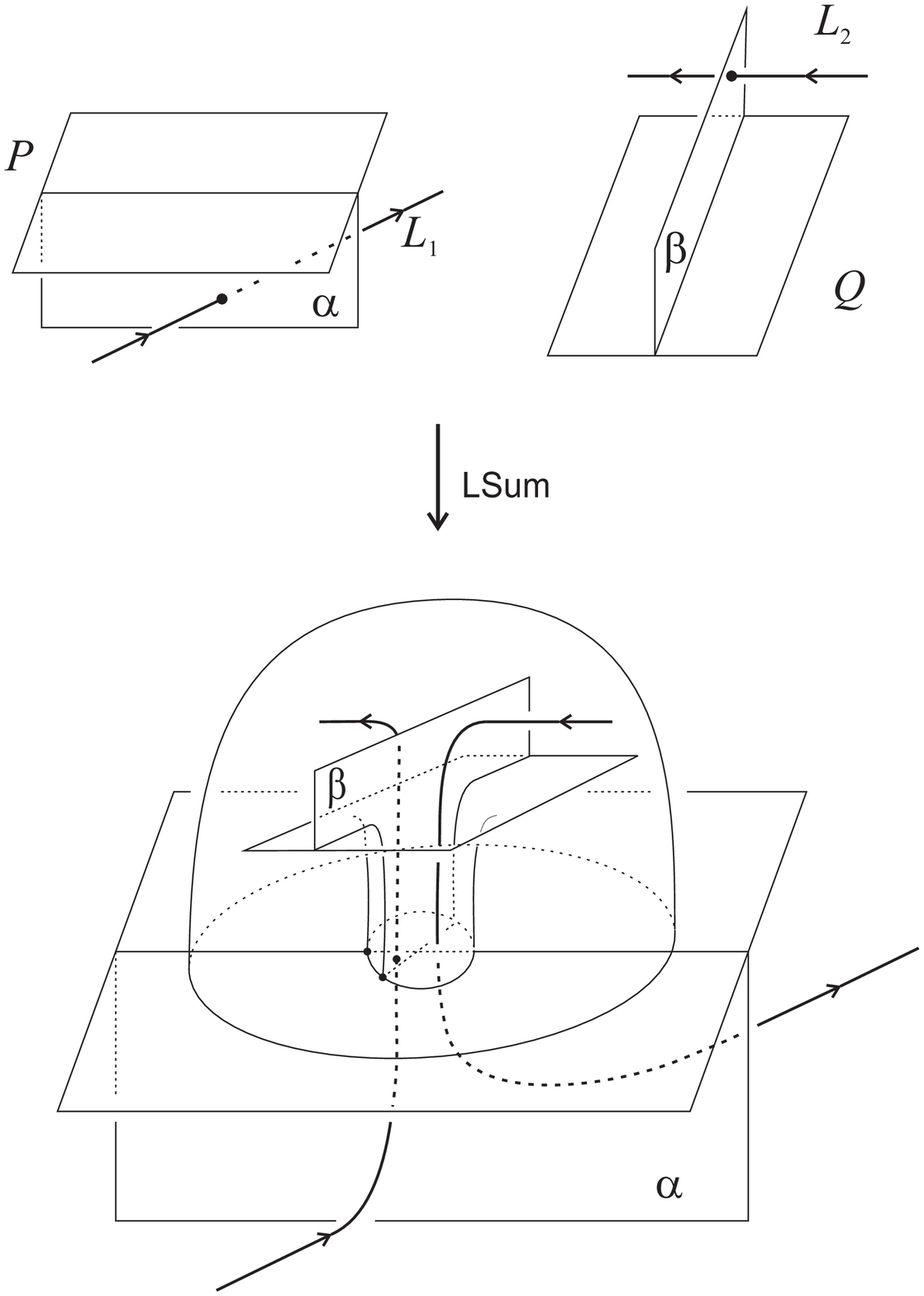}
    \mycap{Construction of $\mbox{{\tt LSum}}_{\alpha,\beta}(P,Q)$ }
    \label{2-conn-sum-spine:fig}
    \end{center}
    \end{figure}

The proof of this lemma is given by
Fig.~\ref{2-conn-sum-spine:fig}, where we assume $s_1=3$ and
$s_2=1$.

\section{State-sums with good initial data}\label{good:data:section}

In this section we describe the behaviour of $\calI_{\xi}$ under
the operations $\mbox{\tt Bub}$, $\mbox{\tt Join}$, $\mbox{\tt
OSum}$, and $\mbox{\tt LSum}$. We then apply these results to
discuss how the invariants behave under ordinary and
link-connected sum of pairs.

\paragraph{Some notation}The following
notation will be handy. If $\alpha$ is a region of a special
polyhedron $P$ and $i$ is a color, we denote by
$\Adm_{(\alpha,i)}(P)$ the set of all admissible colorings of $P$
where $\alpha$ is colored by $i$. Furthermore, if $\xi_1$ and
$\xi_2$ are colorings of disjoint lists of regions of $P$, we
denote their union simply by $\xi_1,\xi_2$. Let now $P$ be an
o-spine of a pair $(M,L)$, let $\alpha_1$, $\ldots$, $\alpha_n$ be
a list of regions of $P$ (that may or may not be marked), and let
$i_1,\ldots,i_n$ be colors. Then for any coloring $\xi$ of the
components of $L$ that are not dual to any of $\alpha_i$, we have
$$\calI_{\xi,(\alpha_1,i_1),\ldots,(\alpha_n,i_n)}(P)=
w^{-m}\sum\limits_{\eta\in\Adm_{\xi,(\alpha_1,i_1),\ldots,(\alpha_n,i_n)}(P)}\omega(\eta).$$
where $m$ is the number of complementary balls having empty
intersection with $L$.

\paragraph{Some properties of good initial data} The first
property that we prove does not in fact require Conditions~$(*)$
--- $(*\!*\!**)$.

\begin{figure}
    \begin{center}
    \includegraphics[scale=0.4]{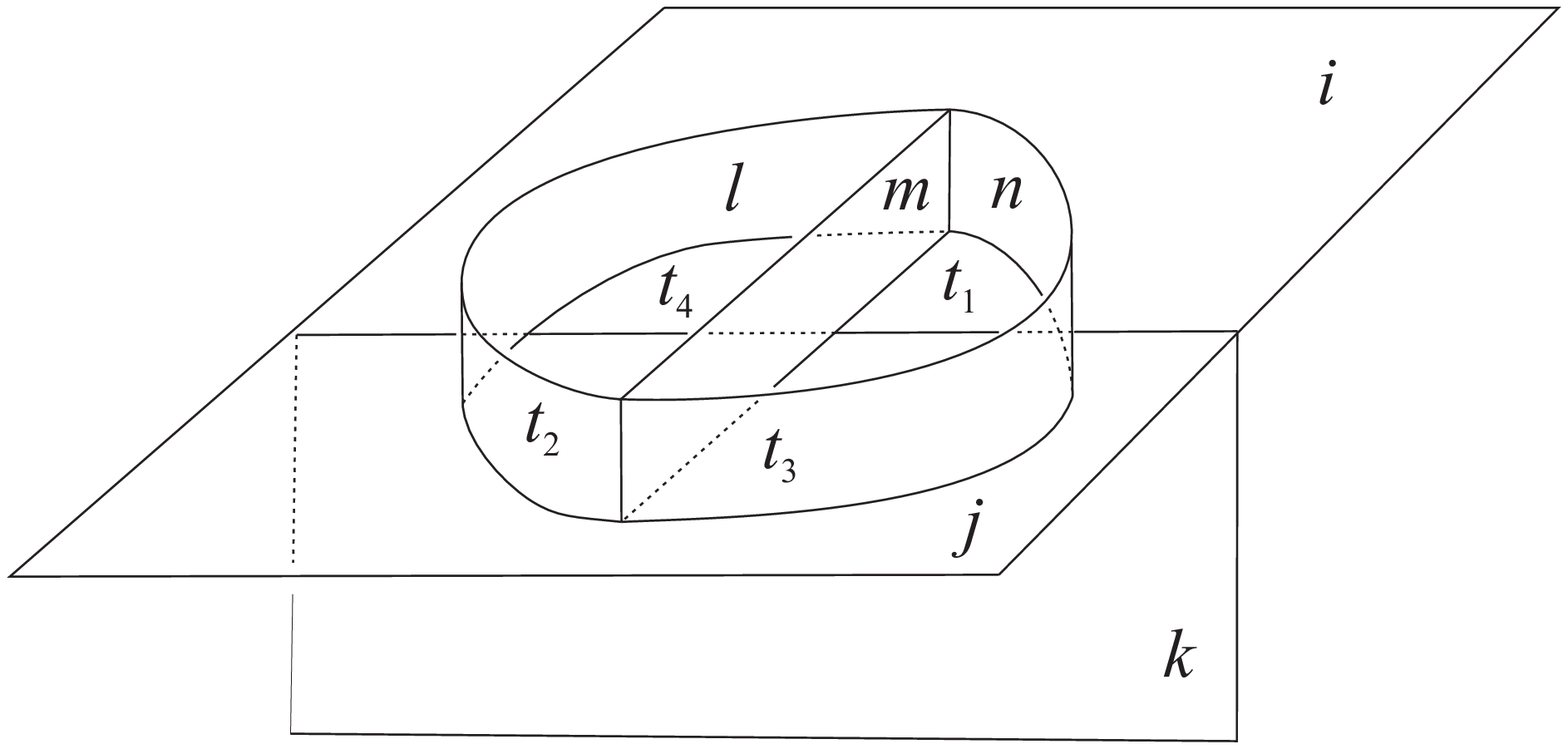}
    \mycap{The polyhedron $P'$}
    \label{strong-irreducibility:fig}
    \end{center}
    \end{figure}

\begin{lemma}\label{strong-irr-wedge:lem}
For any choice of a strongly irreducible datum and for any two
admissible triples $(i,j,k)$ and $(l,m,n)$ there exist
$t_1,t_2,t_3,t_4\in I$ such that the coloring of the polyhedron $P'$ shown in
Fig.~\ref{strong-irreducibility:fig} is admissible.
\end{lemma}

\begin{proof}
Denote by $P''$ the polyhedron obtained from $P'$ by the inverse
$\texttt{MP}$-move along the region colored by $t_1$. By Condition~(iii) of
strong irreducibility the desired $t_1,t_2,t_3,t_4$ exist if and
only if there exist $t_2,t_3,t_4$ such that the resulting coloring
on $P''$ is admissible. Notice that the neighborhoods in $P''$ of
the regions colored by $t_2,t_3,t_4$ have the form shown in
Fig.~\ref{neighbourhoods-of-cells:fig}.
\begin{figure}
    \begin{center}
    \includegraphics[scale=0.6]{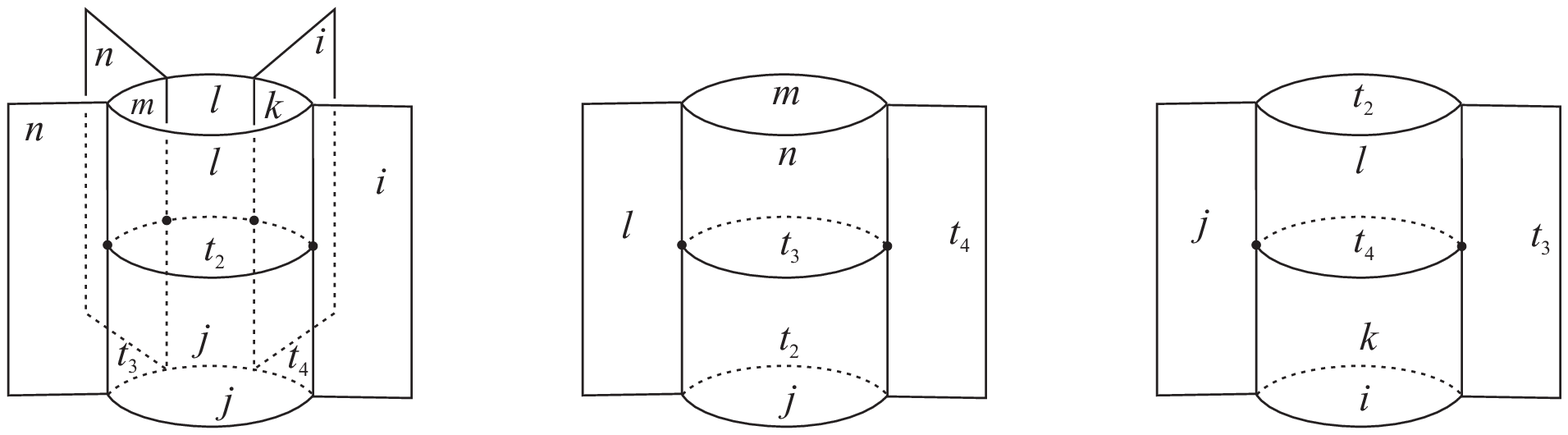}
    \mycap{Neighbourhoods of regions in $P''$}
    \label{neighbourhoods-of-cells:fig}
    \end{center}
    \end{figure}
By Condition~(i) of strong irreducibility we can find $t_2$ such
that the triple $(l,j,t_2)$ is admissible. Since $(l,m,n)$ is
admissible by assumption, by Condition~(ii) of strong irreducibility
we can find a $t_3$ such that both of the triples $(j,n,t_3)$ and
$(m,t_2,t_3)$ are admissible. Similarly, since $(j,i,k)$ and
$(l,j,t_2)$ are admissible, we can find $t_4$ such that $(i,l,t_4)$
and $(k,t_2,t_4)$ are both admissible, whence the conclusion.
\end{proof}

For the rest of this section we now fix some good initial datum.

\begin{lemma}\label{five-eq:lem}
Let $(i,i',i'')$ and $(j,j',j'')$ be admissible triples of colors.
Then
\begin{multline*}
\sum\limits_{k,l,m,n}w_kw_lw_mw_n
 \left|\begin{array}{ccc}i'' & j' & l \\ j & k & j''\end{array}\right|\cdot
 \left|\begin{array}{ccc}i'' & j' & l \\ m & i & i'\end{array}\right| \\
 \times\left|\begin{array}{ccc}i & i'' & i' \\ j'' & n & k\end{array}\right|\cdot
 \left|\begin{array}{ccc}m & j' & i' \\ j'' & n & j\end{array}\right|\cdot
 \left|\begin{array}{ccc}i & l & m \\ j & n & k\end{array}\right|=1
\end{multline*}
with sum taken over all $k,l,m,n$ such that the involved symbols are
defined.
\end{lemma}

\begin{proof}
Collecting together all the summands containing $n$ and applying~$(**)$, we conclude that
that the expression we want to evaluate is equal to
\begin{eqnarray*}
 & & \sum\limits_{k,l,m}w_kw_lw_m
 \left|\begin{array}{ccc}i'' & j' & l \\ j & k & j''\end{array}\right|^2\cdot
 \left|\begin{array}{ccc}i'' & j' & l \\ m & i & i'\end{array}\right|^2\\
   &=& \sum\limits_{k,l}w_kw_l
 \left|\begin{array}{ccc}i'' & j' & l \\ j & k & j''\end{array}\right|^2
 \cdot\left( \sum\limits_m w_m
 \left|\begin{array}{ccc}i'' & j' & l \\ m & i & i'\end{array}\right|^2\right) \\
   &=& \frac{1}{w_{i''}}\sum\limits_{k,l}w_kw_l
 \left|\begin{array}{ccc}i'' & j' & l \\ j & k & j''\end{array}\right|^2
 =\sum\limits_k\frac{w_k}{w_{i''}w_{j''}}
\end{eqnarray*}
where we have used~$(*)$ twice and the last sum is taken over all
$k$ such that $(i'',j'',k)$ is admissible. By Condition~$(*\!*\!**)$
we get the conclusion.
\end{proof}

\begin{figure}
    \begin{center}
    \includegraphics[scale=0.5]{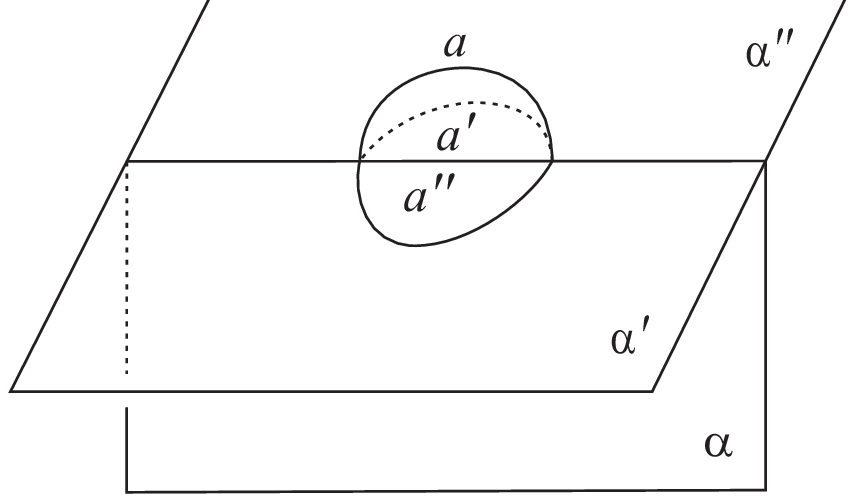}
    \mycap{Notation for the regions of a bubbled polyhedron}
    \label{bubble_2:fig}
    \end{center}
    \end{figure}

\paragraph{State-sums on bubbled polyhedra} Here we
investigate the behaviour of state-sums with good initial data
under the bubble move. The lemmas below use the notation of
Fig.~\ref{bubble_2:fig}.

\begin{lemma}\label{bubble:1:lem}
Let $P$ be an o-spine, and let $\alpha$ be a region of $P$. Let
$\xi$ be a coloring of all marked regions of $P$ except $\alpha$,
$\alpha'$, $\alpha''$ if any of them is marked. Then for any
$i,i',i'',j,k\in I$ such that the triples $(i,i',i'')$ and $(i,j,k)$
are admissible, we have the equality
$$\mathcal{I}_{\xi,(\alpha,i),(\alpha',i'),(\alpha'',i''),(a',j),(a'',k)}\big(\mbox{{\tt Bub}}_{\alpha}(P)\big)=
w^{-1}\frac{w_jw_k}{w_i}\mathcal{I}_{\xi,(\alpha,i),(\alpha',i'),(\alpha'',i'')}(P).$$
\end{lemma}

\begin{proof}
Notice that any coloring in
$\Adm_{\xi,(\alpha,i),(\alpha',i'),(\alpha'',i''),(a',j),(a'',k)}\big(\mbox{{\tt
Bub}}_{\alpha}(P)\big)$ naturally induces a coloring in
$\Adm_{\xi,(\alpha,i),(\alpha',i'),(\alpha'',i'')}(P)$, and the
resulting map
$$\Adm_{\xi,(\alpha,i),(\alpha',i'),(\alpha'',i''),(a',j),(a'',k)}\big(\mbox{{\tt Bub}}_{\alpha}(P)\big)
\rightarrow\Adm_{\xi,(\alpha,i),(\alpha',i'),(\alpha'',i'')}(P)$$ is
surjective (this follows from Condition~(ii) of strong
irreducibility). Notice further that the complement of $\mbox{{\tt Bub}}_{\alpha}(P)$
contains one more empty ball with respect to $P$.
Let $l$ denote the (varying) color of $a$. By the
above remarks we have
\begin{multline*}
\mathcal{I}_{\xi,(\alpha,i),(\alpha',i'),(\alpha'',i''),(a',j),(a'',k)}\big(\mbox{{\tt Bub}}_{\alpha}(P)\big)\\
=w^{-1}\cdot w_jw_k\cdot\left(\sum\limits_lw_l\left|\begin{array}{ccc}i & i' & i'' \\ l &
j &
k\end{array}\right|^2\right)\cdot\mathcal{I}_{\xi,(\alpha,i),(\alpha',i'),(\alpha'',i'')}(P).
\end{multline*}
By Condition~$(*)$ the sum between brackets in this expression
evaluates to $w_i^{-1}$, whence the conclusion.
\end{proof}

\begin{cor}\label{bubble:1:cor}
Under the assumptions of the lemma, we have
$$\mathcal{I}_{\xi,(\alpha,i),(a',j),(a'',k)}\big(\mbox{{\tt Bub}}_{\alpha}(P)\big)=
w^{-1}\frac{w_jw_k}{w_i}\mathcal{I}_{\xi,(\alpha,i)}(P).$$
\end{cor}

\begin{lemma}\label{bubble:2:lem}
Let $P$ be an o-spine, and let $\alpha$ be a region of $P$. Then for
any choice of a coloring $\xi$ of the marked regions of $P$ we have
$$\mathcal{I}_{\xi}\big(\mbox{{\tt Bub}}_{\alpha}(P)\big)=
\mathcal{I}_{\xi}(P).$$
\end{lemma}

\begin{proof}
If $\alpha$ is not marked then by Corollary~\ref{bubble:1:cor} we have that
$$\mathcal{I}_{\xi}\big(\mbox{{\tt Bub}}_{\alpha}(P)\big)=
w^{-1}\cdot\sum\limits_{i,j,k}\frac{w_jw_k}{w_i}\mathcal{I}_{\xi,(\alpha,i)}(P)$$
with sum taken over all $i,j,k$ such that the triple $(i,j,k)$ is
admissible, and the conclusion follows from Condition~$(***)$
together with the obvious equality
$\mathcal{I}_{\xi}(P)=\sum\limits_i\mathcal{I}_{\xi,(\alpha,i)}(P)$.
If $\alpha$ is marked and $\xi'$ is the restriction of $\xi$ to
all the marked regions of $P$ other than $\alpha$ then again by
Corollary~\ref{bubble:1:cor} we have that
$$\mathcal{I}_{\xi}\big(\mbox{{\tt Bub}}_{\alpha}(P)\big)=
w^{-1}\cdot\sum\limits_{j,k}\frac{w_jw_k}{w_i}\mathcal{I}_{\xi',(\alpha,i)}(P)=
w^{-1}\cdot\sum\limits_{j,k}\frac{w_jw_k}{w_i}\mathcal{I}_{\xi}(P),$$
and the conclusion again follows from Condition~$(*\!*\!*)$.
\end{proof}

Using this lemma we can now spell out the proof of the fact that colored Turaev-Viro invariants
for links are well-defined under our assumption that the initial data are good:

\dimo{invt:well-defined:prop} It follows from \cite[Theorem~2.1,
Corollary~2.2]{Amendola} that there exists a sequence of bubble
moves and of admissible $\texttt{Lune}$- and $\texttt{MP}$-moves and their
inverses taking
$P$ to $P'$. The conclusion easily follows, because
Lemma~\ref{bubble:2:lem} implies invariance of the colored
state-sum under bubble moves, whereas, just as in~\cite{TV,Ma},
$(*)$ implies invariance under $\texttt{Lune}$ and its inverse,
and $(**)$ implies invariance under $\texttt{MP}$ and its inverse. \finedimo

\paragraph{State-sums on joins} We now turn to the
behaviour of the invariants under joins of o-spines. We will in
fact consider a somewhat specific situation, which will suffice
for our purposes. Namely, given a germ $\alpha$ of a region along
an edge we consider the three germs of the complementary balls
whose closures have non-empty intersection with $\alpha$. The
closures of two of these germs contain the whole of $\alpha$,
while the closure of the third one intersects it only along a
segment. We will call the ball containing this third germ the
\emph{complementary ball opposite to the germ $\alpha$}.

\begin{prop}\label{wedge:prop}
Let $P$ and $Q$ be o-spines, let $\alpha$ be a germ of a region of
$P$ such that the complementary ball opposite to $\alpha$ is
empty, and let $\beta$ be a germ of a region of $Q$. Let us fix a
coloring $\xi_1$ of the marked regions of $P$ and a coloring
$\xi_2$ of the marked regions of $Q$. Then
$$\mathcal{I}_{\xi_1,\xi_2}\big(\mbox{{\tt
Join}}\big((P,\alpha),(Q,\beta)\big)\big)=
w\cdot \mathcal{I}_{\xi_1}(P)\cdot\mathcal{I}_{\xi_2}(Q).$$
\end{prop}

\begin{proof}
Let $m_1$ be the number of empty complementary balls of $P$, and let $m_2$ be
the number of empty complementary balls of $Q$. By the assumptions of the lemma,
the number of empty complementary balls of $\mbox{{\tt
Join}}\big((P,\alpha),(Q,\beta)\big)$ is equal to $m_1+m_2-1$. Furthermore,
since an admissible coloring of $\mbox{{\tt
Join}}\big((P,\alpha),(Q,\beta)\big)$ induces admissible colorings
on $P$ and on $Q$, there is a naturally defined map
$$\Adm_{\xi_1,\xi_2}\big(\mbox{{\tt
Join}}\big((P,\alpha),(Q,\beta)\big)\big)
\rightarrow\Adm_{\xi_1}(P)\times\Adm_{\xi_2}(Q),$$ and by
Lemma~\ref{strong-irr-wedge:lem} this map is surjective.
\begin{figure}
    \begin{center}
    \includegraphics[scale=0.4]{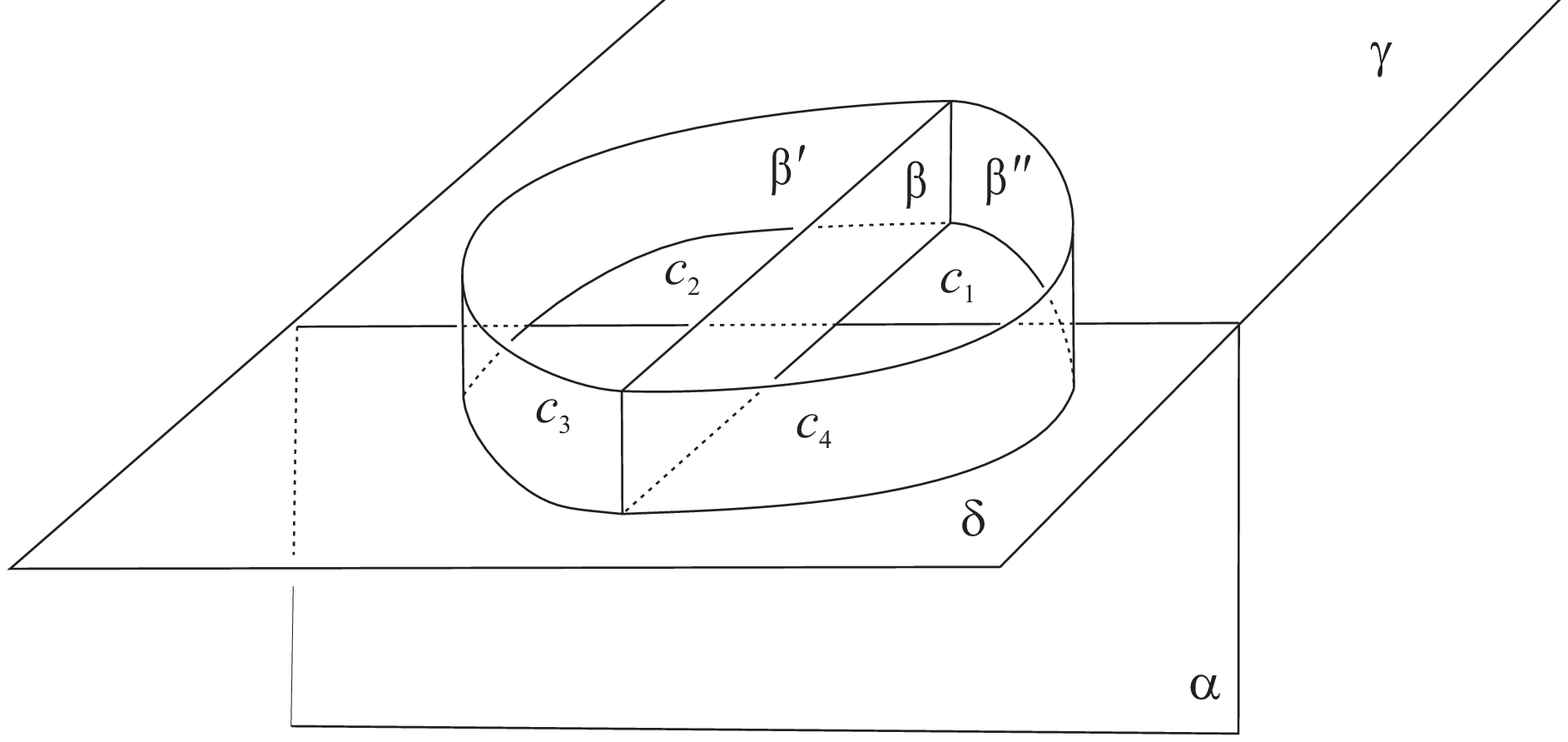}
    \mycap{Notation for the regions of a join}
    \label{wedge-spine-cell-notation:fig}
    \end{center}
    \end{figure}
Using the notation of Fig.~\ref{wedge-spine-cell-notation:fig} and
agreeing that $t_i$ will denote the color of $c_i$ for
$i=1,\ldots,4$, we see that
\begin{gather*}
\mathcal{I}_{\xi_1,\xi_2}\left(\mbox{{\tt
Join}}\left((P,\alpha),(Q,\beta)\right)\right)=
w\cdot w^{-m_1-m_2}\sum\limits_{(\eta_1,\eta_2)}\omega(\eta_1)\cdot\omega(\eta_2)\\
\times\Bigg(\sum\limits_{t_1,t_2,t_3,t_4}w_{t_1}w_{t_2}w_{t_3}w_{t_4}
 \left|\begin{array}{ccc}i'' & j' & t_2 \\ j & t_1 & j''\end{array}\right|\cdot
 \left|\begin{array}{ccc}i'' & j' & t_2 \\ t_3 & i & i'\end{array}\right|\\
\cdot \left|\begin{array}{ccc}i & i'' & i' \\ j'' & t_4 &
t_1\end{array}\right|\cdot
 \left|\begin{array}{ccc}t_3 & j' & i' \\ j'' & t_4 & j\end{array}\right|\cdot
 \left|\begin{array}{ccc}i & t_2 & t_3 \\ j & t_4 & t_1\end{array}\right| \Bigg),
\end{gather*}
where $(\eta_1,\eta_2)$ runs over the set
$\Adm_{\xi_1}(P)\times\Adm_{\xi_2}(Q)$, the
colors $t_i$ are such that the involved symbols are defined, and
we write for brevity $i=\eta_1(\alpha)$, $i'=\eta_1(\delta)$,
$i''=\eta_1(\gamma)$, $j=\eta_2(\beta)$, $j'=\eta_2(\beta')$,
$j''=\eta_2(\beta'')$. By Lemma~\ref{five-eq:lem} the sum between
brackets in this expression always evaluates to~$1$, whence the
conclusion.
\end{proof}

Recall that the set of marked regions of the polyhedron $\mbox{{\tt
LSum}}_{\alpha,\beta}(P,Q)$ consists of the set of marked regions of
$P$ minus $\alpha$, the set of marked regions of $Q$ minus $\beta$,
and a new region $c_{s_1}$, the definition of which was given
immediately before Lemma~\ref{construct:2-con:spine:lem}.

\begin{prop}\label{wedge-sum:2:prop}
Let $P$ and $Q$ be o-spines, let $\alpha$ be a marked region of
$P$, and let $\beta$ be a marked region of $Q$. Let $\xi_1$ be a
coloring of all the marked regions of $P$ other than $\alpha$, and
let $\xi_2$ be a coloring of all the marked regions of $Q$ other
than $\beta$. Let $k$ be a color. Then
$$\calI_{\xi_1,\xi_2,(c_{s_1},k)}\big(\mbox{{\tt LSum}}_{\alpha,\beta}(P,Q)\big)=
\sum\limits_{i,j}\frac{w_k}{w_iw_j}\calI_{\xi_1,(\alpha,i)}(P)\cdot
\calI_{\xi_2,(\beta,j)}(Q)$$ with sum taken over all colors $i,j$
such that the triple $(i,j,k)$ is admissible.
\end{prop}

\begin{proof}
Without loss of generality we may assume that $s_1=1$ in the
notation of Fig.~\ref{join-cell-notation:fig}. Then $s_2=3$.
Recall that $\mbox{{\tt LSum}}_{\alpha,\beta}(P,Q)$ is obtained
from a certain $\mbox{{\tt Join}}\big((\mbox{{\tt
Bub}}_{\alpha}(P),\overline{\alpha}),(Q,\beta)\big)$ by
removing
the region $c_3$. As a consequence of this removal,
the regions $c_2$, $c_4$, and $a''$
get absorbed into $\alpha$, $\beta$, and $\beta''$, respectively. Thus, an
admissible coloring of $\mbox{{\tt LSum}}_{\alpha,\beta}(P,Q)$
induces admissible colorings on $\mbox{{\tt Bub}}_{\alpha}(P)$ and
on $Q$, and the resulting map of colorings
\begin{multline*}
\Adm_{\xi_1,\xi_2,(c_1,k),(\alpha,i),(a',i'),(\beta,j),(\beta',j'),(\beta'',j'')}
\big(\mbox{{\tt LSum}}_{\alpha,\beta}(P,Q)\big)\\
\longrightarrow\Adm_{\xi_1,(\alpha,i),(a',i'),(b,j')}\big(\mbox{{\tt
Bub}}_{\alpha}(P)\big)\times
\Adm_{\xi_2,(\beta,j),(\beta',j'),(\beta'',j'')}(Q)
\end{multline*}
is onto for any $i$, $j$, $i'$, $j'$, $j''$ such that the triples
$(i,j,k)$, $(i',j',k)$, and $(i,i',j'')$ are admissible.

Denoting the varying color of $\alpha=c_2$ by $i$, the color of $a'$
by $i'$, the color of $\beta=c_4$ by $j$, the color of $\beta'$ by
$j'$, and the color of $\beta''=a''$ by $j''$, from the construction
of $\mbox{{\tt LSum}}_{\alpha,\beta}(P,Q)$ we conclude that
\begin{multline*}
\calI_{\xi_1,\xi_2,(c_1,k)}\big(\mbox{{\tt LSum}}_{\alpha,\beta}(P,Q)\big)
=\sum\limits_{i,j}\sum\limits_{j',j''} \sum\limits_{i'}
w\cdot\calI_{\xi_1,(\alpha,i),(a',i'),(a'',j'')}\big(\mbox{{\tt
Bub}}_{\alpha}(P)\big)\\
\times\calI_{\xi_2,(\beta,j),(\beta',j'),(\beta'',j'')}(Q)\cdot
w_kw_{j''}^{-1}
 \left|\begin{array}{ccc}k & i & j \\ j'' & j' & i'\end{array}\right|^2
\end{multline*}
with sum taken over all indices such that the symbol is defined. By
Corollary~\ref{bubble:1:cor} this is equal to
$$\sum\limits_{i,j}\calI_{\xi_1,(\alpha,i)}(P)\cdot
\sum\limits_{j',j''}
\calI_{\xi_2,(\beta,j),(\beta',j'),(\beta'',j'')}(Q)\cdot
w_kw_i^{-1}\cdot\left(\sum\limits_{i'} w_{i'}
 \left|\begin{array}{ccc}k & i & j \\
j'' & j' & i'\end{array}\right|^2\right).$$ Using Condition~$(*)$, we
get
\begin{multline*}
\sum\limits_{i,j}\calI_{\xi_1,(\alpha,i)}(P)
w_kw_i^{-1}w_j^{-1}\cdot\sum\limits_{j',j''}
\calI_{\xi_2,(\beta,j),(\beta',j'),(\beta'',j'')}(Q)\\=\sum\limits_{i,j}
\frac{w_k}{w_iw_j}\calI_{\xi_1,(\alpha,i)}(P)\cdot\calI_{\xi_2,(\beta,j)}(Q),
\end{multline*}
whence the conclusion.
\end{proof}

\paragraph{Behaviour under connected sum} The above results now easily allow us
to describe the behaviour of the invariants under
connected sum. For the ordinary
version, Proposition~\ref{wedge:prop} and Lemma~\ref{bubble:2:lem}
immediately imply:

\begin{thm}\label{ordinary:connected:thm}
Let $\xi_t$ be a coloring of a link $L_t\subset M_t$ for $t=1,2$.
Then
$$\calI_{\xi_1,\xi_2}(M_1\#_o M_2,L_1\sqcup L_2)=w\cdot
\calI_{\xi_1}(M_1,L_1)\cdot\calI_{\xi_2}(M_2,L_2).$$
\end{thm}

Let us now turn to the link-connected sum. Let $K_t$ be a connected
component of a link $L_t\subset M_t$, let $\xi_t$ be a coloring
of $L_t\setminus K_t$ for $t=1,2$, and let $(M,L)$ be
the link-connected sum of $(M_1,L_1)$ and $(M_2,L_2)$ along $K_1$
and $K_2$. If $\xi$ is one of the colorings $\xi_1$, $\xi_2$, or
$\xi_1,\xi_2$ then there is one component of the relevant link not
colored by $\xi$, and we define the coloring $(\xi;i)$ by giving
it color $i$.
Proposition~\ref{wedge-sum:2:prop} implies:

\begin{thm}\label{2-connected:thm}
We have
$$\frac{1}{w_k}\calI_{(\xi_1,\xi_2;\;k)}(M,L)=\sum\limits_{i,j}\frac{1}{w_i}
\calI_{(\xi_1;\;i)}(M_1,L_1)\cdot\frac{1}{w_j}\calI_{(\xi_2;\;j)}(M_2,L_2)$$ with
sum taken over all $i,j$ such that the triple $(i,j,k)$ is
admissible.
\end{thm}

In the case of knots the above formula has a particularly nice form:
$$\frac{1}{w_k}\calI_k(M_1\# M_2,K_1\# K_2)=
\sum\limits_{i,j}\frac{1}{w_i}\calI_i(M_1,K_1)\cdot\frac{1}{w_j}\calI_j(M_2,K_2).$$

\section{The $\varepsilon$-invariant}\label{eps-invt:0:section}

We now specialize our analysis of colored Turaev-Viro invariants
to the case of the good initial data of~\cite{epsilon}. Another specialization
will be discussed in the next section.

\paragraph{Definition} The so-called $\varepsilon$-invariant
of $3$-manifolds was introduced in \cite{epsilon} (see also
\cite[Section 8.1.3]{Ma}). Let us fix a root $\varepsilon$ of the
equation $x^2=x+1$, and an arbitrary square root
$\varepsilon^{1/2}$ of $\varepsilon$. The corresponding initial datum
$\mathcal{D}_{\varepsilon}$ consists of the following. The ring $\matK$ is $\matC$. The set $I$
is the set $\{0,1\}$. The weights are given by $w_0=1$,
$w_1=\varepsilon$, and the constant $w$ is equal to $\varepsilon+2$.
The admissible triples are those satisfying the so-called
\emph{triangle inequalities}: a triple $(i,j,k)$ is admissible if
and only if
\begin{equation}\label{triang:ineq:eq}
i\leqslant j+k,\;\;j\leqslant k+i,\;\;k\leqslant i+j.
\end{equation}
The values of the $6j$-symbols are given by:
\begin{multline*}
 \left|\begin{array}{ccc}0 & 0 & 0 \\ 0 & 0 & 0\end{array}\right|_{\varepsilon}=1,\;\;
 \left|\begin{array}{ccc}0 & 0 & 0 \\ 1 & 1 & 1\end{array}\right|_{\varepsilon}=
   \varepsilon^{-1/2},\\
 \left|\begin{array}{ccc}0 & 1 & 1 \\ 0 & 1 & 1\end{array}\right|_{\varepsilon}=
   \varepsilon^{-1},\;\;
 \left|\begin{array}{ccc}0 & 1 & 1 \\ 1 & 1 & 1\end{array}\right|_{\varepsilon}=
   \varepsilon^{-1},\;\;
 \left|\begin{array}{ccc}1 & 1 & 1 \\ 1 & 1 & 1\end{array}\right|_{\varepsilon}=
   -\varepsilon^{-2}.
\end{multline*}

\begin{prop}\label{eps-good:lem}
The datum just described is good.
\end{prop}

\begin{proof}
Strong irreducibility and Condition~$(*\!*\!*\,*)$ can be checked directly.
Condition~$(**)$ was established in~\cite[Section 8.1.3]{Ma}.

Let us prove that the datum satisfies~$(*\!*\!*)$. If $k=0$
then the sum on the left-hand side of~$(*\!*\!*)$ is equal to
$w_0^2+w_1^2=\varepsilon+2=(\varepsilon+2)w_0$. If $k=1$ then the
sum is equal to
$w_0w_1+w_1w_0+w_1^2=(\varepsilon+2)\varepsilon=w\cdot w_1$,
whence the conclusion.

To prove orthogonality, we consider various cases. If $i=i'=0$,
we must evaluate
$$\sum\limits_mw_m\left|\begin{array}{ccc}0 & j & j \\ m & l &
l\end{array}\right|_{\varepsilon}^2,$$ where up to symmetry we can
assume that $(j,l)$ is $(0,0)$, $(1,0)$, or $(1,1)$, and a simple calculation
shows that the sum is always $1$, whence the desired conclusion.

Suppose now that $i=i'=1$. Then up to symmetry $(j,k)$ is either
$(1,0)$ or $(1,1)$. If $(j,k)=(1,0)$ then $m=l$, so we should
evaluate the expression $$w_l\left|\begin{array}{ccc}1 & 1 & 0 \\
l & l & n\end{array}\right|_{\varepsilon}^2$$ for $l,n$
not simultaneously zero. In all cases we get
$\varepsilon^{-1}$, \emph{i.e.} $w_1^{-1}$. Suppose now that
$(j,k)=(1,1)$. If either $l$ or $n$ is zero, the value of the
expression is
$$w_1\left|\begin{array}{ccc}1 & 1 & 1 \\ 1 & 1 &
0\end{array}\right|_{\varepsilon}^2=\varepsilon^{-1}=w_1^{-1}$$ while
if also $l=n=1$, it is
$$w_0\left|\begin{array}{ccc}1 & 1 & 1 \\ 0 & 1 &
1\end{array}\right|_{\varepsilon}^2+
 w_1\left|\begin{array}{ccc}1 & 1 & 1 \\ 1 & 1 & 1\end{array}\right|_{\varepsilon}^2,$$
which gives $\varepsilon^{-1}$ as well.

Finally, suppose that $i=0$ and $i'=1$. Notice that $(0,j,j)$ and
$(1,j,j)$ are both admissible only if $j=1$. For the same reason
$l=1$, and the left-hand side of~$(*)$ is
$$\left|\begin{array}{ccc}0 & 1 & 1 \\ 0 & 1 &
1\end{array}\right|_{\varepsilon}\cdot\left|\begin{array}{ccc}1 &
1 & 1 \\ 0 & 1 & 1\end{array}\right|_{\varepsilon}+
w_1\left|\begin{array}{ccc}0 & 1 & 1 \\ 1 & 1 &
1\end{array}\right|_{\varepsilon}\cdot\left|\begin{array}{ccc}1 &
1 & 1 \\ 1 & 1 & 1\end{array}\right|_{\varepsilon},$$ which is
indeed equal to zero.
\end{proof}

Let us denote by $\calI_{*}^{\varepsilon}$ the colored Turaev-Viro
invariants of links obtained with respect to the datum just
described.

\paragraph{A refinement of the formula for link-connected sums}
Theorem~\ref{2-connected:thm} yields the following result, that will lead
to important applications in Section~\ref{other-invts:section}.

\begin{thm}\label{eps:2-connected:thm}
Let $(M,L)$ be the link-connected sum of $(M_1,L_1)$ and $(M_2,L_2)$
along knots $K_1\subset L_1$ and $K_2\subset L_2$. Let $L_i'$ be
$L_i\setminus K_i$. Then
\begin{eqnarray*}
\calI_0^{\varepsilon}(M,L)& = & \calI_0^{\varepsilon}(M_1,L_1)\cdot
\calI_0^{\varepsilon}(M_2,L_2) \\
& + & \varepsilon^{-2}\cdot
\left((\varepsilon+2)\calI_0^{\varepsilon}(M_1,L_1')-\calI_0^{\varepsilon}(M_1,L_1)\right)\\
& & \ \,\quad \times
\left((\varepsilon+2)\calI_0^{\varepsilon}(M_2,L_2')-\calI_0^{\varepsilon}(M_2,L_2)\right).
\end{eqnarray*}
\end{thm}

\begin{proof}
By Theorem~\ref{2-connected:thm} we have
$$\calI_0^{\varepsilon}(M,L)=
\calI_0^{\varepsilon}(M_1,L_1)\cdot\calI_0^{\varepsilon}(M_2,L_2)+
\varepsilon^{-2}\calI_{(0;K_1,1)}^{\varepsilon}(M_1,L_1)\cdot
\calI_{(0;K_2,1)}^{\varepsilon}(M_2,L_2).$$ Moreover, by the definition of
$\calI_{*}^{\varepsilon}$ we obtain
$$\calI_{(0;K_i,0)}^{\varepsilon}(M_i,L_i)+
\calI_{(0;K_i,1)}^{\varepsilon}(M_i,L_i)=w\calI_0^{\varepsilon}(M_i,L_i')
=(\varepsilon+2)\calI_0^{\varepsilon}(M_i,L_i'),$$
whence the conclusion.
\end{proof}

\paragraph{The $\varepsilon$-invariant of unlinks}
Theorem~\ref{ordinary:connected:thm} gives the following:

\begin{cor}\label{eps:unlink:cor}
Let $U_n$ be the unlink with $n\geqslant 0$ components in $\matS^3$. Then
$$\calI_0^{\varepsilon}(U_n)=(\varepsilon+2)^{n-1}.$$
\end{cor}

\begin{proof}
By definition
$\calI_0^{\varepsilon}(\emptyset)=t(\matS^3)=(\varepsilon+2)^{-1}$.
Furthermore, it has already been mentioned that Bing's house with the region $c$
marked (see Fig.~\ref{Abalone-Bing:fig}-right) is an o-spine of the
unknot. If this region is colored by $0$, the only admissible
coloring is when the other two regions are also colored by $0$. Hence
$\calI_0^{\varepsilon}(U_1)=1$. The desired conclusion now follows from
Theorem~\ref{ordinary:connected:thm}.
\end{proof}

\begin{cor}\label{eps:non-zero:color:cor}
Let $\xi$ be a coloring of $U_n$ with at least one color $1$. Then
$$\calI_{\xi}^{\varepsilon}(U_n)=0.$$
\end{cor}

\begin{proof}
Using the abalone, it is easy to calculate directly that
$\calI_1^{\varepsilon}(U_1)=0$. The conclusion now follows from
Theorem~\ref{ordinary:connected:thm}.
\end{proof}

\paragraph{$\calI_0^{\varepsilon}$ is the $\varepsilon$-invariant of the complement}
Corollaries~\ref{eps:unlink:cor}
and~\ref{eps:non-zero:color:cor} suggest that the invariant
$\calI_0^{\varepsilon}$ is well-normalized and so is perhaps the
most natural one to consider. On the other hand, as we prove in
this paragraph, it has the limitation of being determined by the link complement.

\begin{lemma}\label{eps:symbol-with-0:lem}
For any admissible triples $(i,j,k)$ and $(i,j,m)$ we have
$$\left|\begin{array}{ccc}0 & i & i \\ k & j &
j\end{array}\right|_{\varepsilon}\cdot
 \left|\begin{array}{ccc}0 & i & i \\ m & j & j\end{array}\right|_{\varepsilon}
 =(w_iw_j)^{-1}.$$
\end{lemma}

\begin{proof}
The proof is by inspection of the possibilities.
\end{proof}

\begin{lemma}\label{eps-0-puncture:lem}
Let $P$ be an o-spine, and let $c$ be a length-$2$ region of $P$
with embedded closure. Suppose that for either edge $e$ of $c$
the two germs of region incident to $e$ and not contained
in $c$ are contained in two distinct regions. Suppose also that
$c$ has a common edge with at least three distinct regions. Let $P'$ be the marked
polyhedron obtained from $P$ by removing $c$, and let $\xi'$ be a coloring
of the marked regions of $P'$. Then there is a bijection
$\imath:\Adm_{\xi',(c,0)}(P)\rightarrow\Adm_{\xi'}(P')$ such that
for any $\eta\in\Adm_{\xi',(c,0)}(P)$ we have
$$\omega(\eta)=\omega(\imath(\eta)).$$
\end{lemma}

\begin{proof}
The assumptions of the lemma imply that $P'$ is a special polyhedron.
We observe that each region of $P'$ either corresponds to a unique
region of $P$ or to the union of several (at most three) regions of
$P$. Notice that in the latter case for each region $c_1$ in the
union there is a region $c_2$ also in the union such the triple $c_1,c_2,c$ has a
common edge. Let now $\eta\in\Adm_{\xi',(c,0)}(P)$. Since $(0,0,1)$
is not admissible, if two regions of $P$ are contained in one
region of $P'$ then they have the same color. So $\eta$ defines an admissible coloring
$\imath(\eta)$ of $P'$, and it is easy to see that this correspondence is bijective.

Since $c$ has length 2, its admissibly colored neighbourhood has the form shown in
Fig.~\ref{length-2-cell:fig}.
\begin{figure}
    \begin{center}
    \includegraphics[scale=0.5]{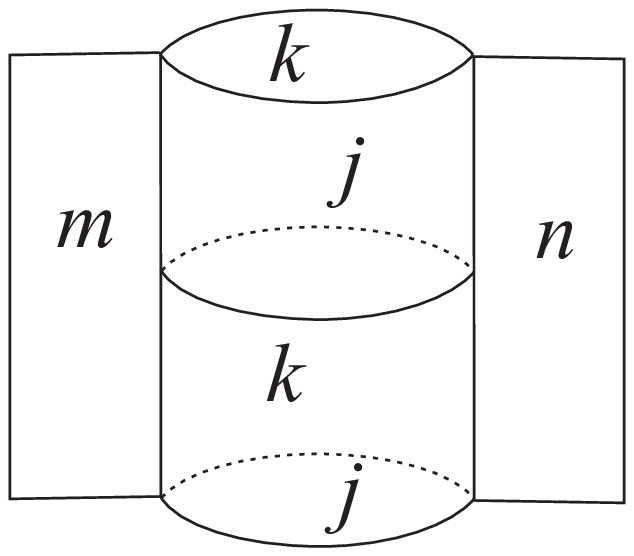}
    \mycap{A colored neighbourhood of a length-2 region with embedded closure}
    \label{length-2-cell:fig}
    \end{center}
    \end{figure}
Notice that, although some of the germs in Fig.~\ref{length-2-cell:fig} may belong
to the same region, the total number of regions in $P'$ is smaller than that of
$P$ by $3$. Hence we see that
$$\omega(\imath(\eta))=w_j^{-1}w_k^{-1}\left(
\left|\begin{array}{ccc}0 & j & j \\ m & k &
k\end{array}\right|_{\varepsilon}\cdot \left|\begin{array}{ccc}0 &
j & j \\ n & k & k\end{array}\right|_{\varepsilon}
\right)^{-1}\omega(\eta),$$ and, as before, we need to ensure that
the factor is always equal to $1$. Since the triples $(j,k,m)$ and
$(j,k,n)$ must be admissible, the conclusion follows from
Lemma~\ref{eps:symbol-with-0:lem}.
\end{proof}

\begin{lemma}\label{good-o-spine:lem}
Every pair $(M,L)$ admits an o-spine $P$ such that each
marked region of $P$
satisfies the conditions of Lemma~\ref{eps-0-puncture:lem}.
\end{lemma}

\begin{proof}
Let us choose an arbitrary o-spine $P'$ of $M$ and a projection
$\Pi$ of $L$ on $P$. Of course we can and will assume that the
projection has only generic intersections with itself and with
$S(P')$ and that it decomposes each region of $P'$ into discs. We
then add a curl (as in the first Reidemeister move) on each
component of $\Pi$ that has no self-intersection. We then proceed
exactly as in the construction of $P_{\Pi}(L)$ for links in
$\matS^3$, digging a tunnel inside $P$ along $\Pi$, adding marked
length-2 regions but omitting the last step (that of
Fig.~\ref{1-move-proj-spine:fig}). It is easy to check that the
resulting polyhedron $P$ is an o-spine of $(M,L)$ satisfying all
the desired conditions.
\end{proof}

We now recall that the $\varepsilon$-invariant~\cite{epsilon} of a $3$-manifold $M$
arising from the initial datum $\mathcal{D}_{\varepsilon}$ is denoted by
$t(M)$.

\begin{thm}\label{eps:link-compl:thm}
Let $L\subset M$ be a link in a closed $3$-manifold $M$, and let $N_M(L)$ denote its
complement. Then
$$\calI_0^{\varepsilon}(M,L)=t(N_M(L)).$$
\end{thm}

\begin{proof}
Let $P$ be an o-spine of $(M,L)$ as in Lemma~\ref{good-o-spine:lem}, let $m$ be
the number of empty balls in the complement $M\setminus P$, and let
$P'$ be the polyhedron obtained from $P$ by removing all the marked components. Then
$P'$ is a special skeleton (see Remark~\ref{invt:for:bry:rem}) of $N_M(L)$, with
$N_M(L)\setminus P'$ containing a collar of the boundary and $m$
balls. By repeated application
of Lemma~\ref{eps-0-puncture:lem} there exists a bijection
$\imath:\Adm_0(P)\rightarrow\Adm(P')$ such that for any $\eta\in\Adm_0(P)$ we have
$\omega(\eta)=\omega(\imath(\eta))$. By definition,
$\calI_0^{\varepsilon}(L)=w^{-m}\sum\limits_{\eta\in\Adm_0(P)}\omega(\eta)$ and
$$t(N_M(L))=w^{-m}\sum\limits_{\eta'\in\Adm(P')}\omega(\eta')=w^{-m}\sum\limits_{\eta\in\Adm_(P)}\omega(\imath(\eta)),$$
whence the conclusion.
\end{proof}

In Section~\ref{invt-not-complement:sect} we will show that the
invariant $\calI_1^{\varepsilon}$ does not have the drawback of
being determined by the link complement.

\begin{figure}
    \begin{center}
    \includegraphics[scale=0.5]{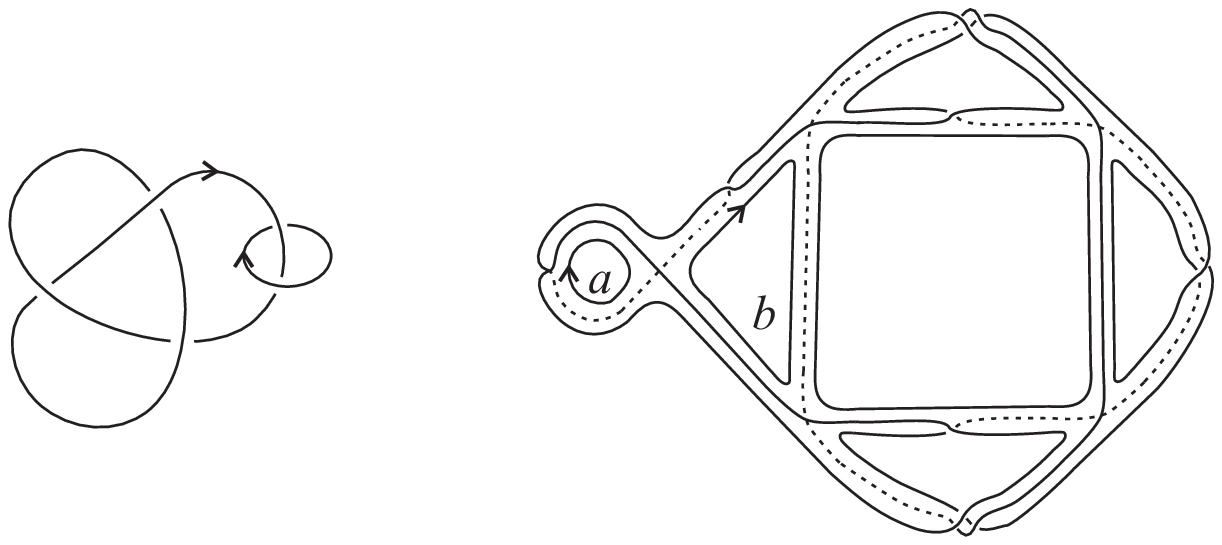}
    \mycap{A low-crossing link $L$ and an o-spine of $L$. Note that
    the trefoil component of $L$ is dual to the region $a$ of $P$ and
    the unknot component of $L$ is dual to the region $b$}
    \label{link-tref-triv:fig}
    \end{center}
    \end{figure}

\paragraph{Examples}
It follows from Corollary~\ref{eps:unlink:cor} that
for the unknot $U_1$ we have
$$\calI_0^{\varepsilon}(U_1)=1.$$ Using the abalone
(shown in Fig.~\ref{Abalone-Bing:fig}-left) with the region $b$
marked to calculate the invariant of the trefoil $3_1$, we
easily see that
$$\calI_0^{\varepsilon}(3_1)=2.$$
We also compute $\calI_0^{\varepsilon}$ for the link $L$ shown
in Fig.~\ref{link-tref-triv:fig}-left.
An o-spine of $L$ is shown in Fig.~\ref{link-tref-triv:fig}-right, and
a direct calculation yields
$$\calI_0^{\varepsilon}(L)=5-\varepsilon.$$

\section{Turaev-Viro invariants}

We will now discuss the Turaev-Viro initial
datum introduced in \cite[Section~7.1]{TV}.

\paragraph{Definition} First of all, we fix an integer $r\geqslant 3$
(called the \emph{order}
of the datum and of the corresponding invariant) and a root of unity
$q$ of degree $2r$ such that $q^2$ is a primitive root of unity of
degree $r$. Recall that the \emph{quantum factorial} $[n]!$ of an
integer $n\geqslant 1$ is defined by the following equalities:
$$[n]=\frac{q^{n}-q^{-n}}{q-q^{-1}},\;\;\;[n]!=[n][n-1]\ldots[1].$$
Also by definition $[0]!=[0]=1$.

Let $\matK$ be $\matC$ and let $I$ be the set $\{0,1/2,1,\ldots,(r-3)/2,(r-2)/2\}$. A triple
$(i,j,k)$ is admissible if it satisfies the triangle
inequalities~(\ref{triang:ineq:eq}) and the following two
conditions:
$$i+j+k\in\mathbb{N}, \mbox{ and~ }i+j+k\leqslant r-2.$$
The color weights and the constant $w$ are given by
$$w_i=(-1)^{2i}[2i+1], \;\;\;\; w=-\frac{2r}{(q-q^{-1})^2}.$$

To define the $6j$-symbols, which are called \emph{quantum
$6j$-symbols}, we need the following additional functions. For every
admissible triple $(i,j,k)$ we set
$$\Delta(i,j,k)=\left(\frac{[i+j-k]![j+k-i]!
[k+i-j]!}{[i+j+k+1]!}\right)^{1/2},$$ where by the square root of a
real number $a$ we mean $a^{1/2}\geqslant 0$ if $a\geqslant 0$ and
$\sqrt{-1}\,|a|^{1/2}$ if $a<0$. For every admissible 6-tuple
$(i,j,k,l,m,n)$ we also set
\begin{multline*}
U(z,(i,j,k,l,m,n))=[z-i-j-k]![z-i-m-n]!\\
\shoveright{\times[z-j-l-n]![z-k-l-m]!}\\
\shoveleft{L(z,(i,j,k,l,m,n))=[i+j+l+m-z]!}\\
\times[j+k+m+n-z]![k+i+n+l-z]!
\end{multline*}
for any $z$ such that all the factorials are defined. Now we define:
\begin{multline*}
\left|\begin{array}{ccc}i & j & k \\ l & m &
n\end{array}\right|=(\sqrt{-1})^{-2(i+j+k+l+m+n)}\Delta(i,j,k)\Delta(i,m,n)
\Delta(j,l,n)\Delta(k,l,m)\\
\times\sum\limits_z(-1)^z
[z+1]!U(z,(i,j,k,l,m,n))^{-1}L(z,(i,j,k,l,m,n))^{-1},
\end{multline*}
where the sum is taken over all $z$ such that
\begin{multline*}
\max\{i+j+k,i+m+n,j+l+n,k+l+m\}\leqslant z,\\
z\leqslant\min\{i+j+l+m,j+k+m+n,k+i+n+l\}.
\end{multline*}

We notice that the first condition of strong irreducibility and
Conditions~$(*)$, $(**)$, $(*\!*\!*)$ were established in the proof of
\cite[Theorem~7.2]{TV}.

\begin{lemma}\label{tv-almost-good:lem}
The initial datum just described is good.
\end{lemma}

\begin{proof}
To check the second condition of strong irreducibility, we need to
show that, given admissible triples $(i,j,k)$, $(i,m,n)$, there
exists $l$ such that $(j,m,l)$ and $(k,n,l)$ are also admissible.
Admissibility of $(j,m,l)$ means:
\begin{itemize}
  \item $|m-j|\leqslant l\leqslant\min\{r-2-(m+j),m+j\}$;
  \item $l\equiv m+j\pmod 1$.
\end{itemize}
Similarly, admissibility of $(k,n,l)$ means
\begin{itemize}
  \item $|n-k|\leqslant l\leqslant\min\{r-2-(k+n),k+n\}$;
  \item $l\equiv k+n\pmod 1$.
\end{itemize}
Without loss of generality, we may and will assume that $m\geqslant
j$ and $m-j\geqslant|n-k|$. We remark that $m+j\equiv k+n\pmod 1$,
and that $m-j\leqslant\min\{r-2-(m+j),m+j\}$, because
$m\leqslant(r-2)/2$, being a member of an admissible triple.
Thus, for the desired $l$ to exist it is sufficient to show that
\begin{itemize}
  \item $m-j\leqslant\min\{r-2-(k+n),k+n\}$.
\end{itemize}
Using admissibility of $(i,j,k)$ and $(i,m,n)$, we see that
$$m-j\leqslant i+n-j\leqslant j+k+n-j=n+k,$$
$$m-j+n+k\leqslant m-j+n+i+j=i+m+n\leqslant r-2.$$
An identical argument shows that there exists $l'$ such that
$(k,m,l')$ and $(j,n,l')$ are admissible, whence the conclusion.

Let us now check the third condition. Notice that
admissibility of the 6-tuples $(i,j,k,j_1,k_1,i_1)$ and $(i,j,k,j_2,k_2,i_2)$
is equivalent to admissibility
of the triples $(i,j,k)$, $(i,i_1,k_1)$, $(i,i_2,k_2)$,
$(j,i_1,j_1)$, $(j,i_2,j_2)$, $(k,j_1,k_1)$, $(k,j_2,k_2)$. Therefore
$$i_1+i_2\equiv 2i+k_1+k_2\equiv k_1+k_2\pmod 1,\;\;
j_1+j_2\equiv 2j+i_1+i_2\equiv i_1+i_2\pmod 1,$$
so for the necessary $l$ to
exist it is sufficient to show that
\begin{multline*}
\max\{|i_1-i_2|,|j_1-j_2|,|k_1-k_2|\}\\
\leqslant\min\{i_1+i_2,j_1+j_2,k_1+k_2,r-2-(i_1+i_2),r-2-(j_1+j_2),r-2-(k_1+k_2)\}.
\end{multline*}
Indeed, without loss of generality we will assume that
$\max\{|i_1-i_2|,|j_1-j_2|,|k_1-k_2|\}=k_1-k_2$. Then we have
$$k_1-k_2\leqslant i+i_1-k_2\leqslant i_2+k_2+i_1-k_2=i_1+i_2,$$
$$k_1-k_2\leqslant k+j_1-k_2\leqslant j_2+k_2+j_1-k_2=j_1+j_2,$$
and of course $k_1-k_2\leqslant k_1+k_2$. Furthermore, we have
$$k_1-k_2+i_1+i_2\leqslant k_1-k_2+i_1+i+k_2=k_1+i_1+i\leqslant r-2,$$
$$k_1-k_2+j_1+j_2\leqslant k_1-k_2+j_1+k+k_2=k_1+j_1+k\leqslant r-2,$$
$$k_1-k_2+k_1+k_2=2k_1\leqslant r-2,$$
whence the conclusion.

To conclude, we need to check Condition~$(*\!*\!*\,*)$. Since by
\cite[Equation~(19)]{TV} we have
$$\frac{1}{w_iw_j}=\left|\begin{array}{ccc}0 & i & i \\ k & j &
j\end{array}\right|^2,$$ rewriting the desired equality as
$\sum\limits_k\frac{w_k}{w_iw_j}=1$, we obtain
$$\sum\limits_k\frac{w_k}{w_iw_j}=\sum\limits_kw_k\left|\begin{array}{ccc}0 & i & i \\ k & j &
j\end{array}\right|^2=w_0^{-1}=1,$$ where the latter two equalities
are respectively by Condition~($*)$ and by direct calculation.
\end{proof}

We denote the colored Turaev-Viro invariants corresponding to a given choice of
$r,q$ by $\calI_*^{r,q}$.

\paragraph{The invariants of unlinks} The following is analogous to
Corollary~\ref{eps:unlink:cor}:

\begin{cor}\label{tv:unlink:cor}
For any $n\geqslant 0$ we have
$$\calI_0^{r,q}(U_n)=\left(-\frac{2r}{(q-q^{-1})^2}\right)^{n-1}.$$
\end{cor}

\begin{proof}
By definition
$\calI_0^{r,q}(\emptyset)=TV(\matS^3)=\left(-\frac{2r}{(q-q^{-1})^2}\right)^{-1}$.
Let us show that $\calI_0^{r,q}(U_1)=1$ for any $r$. Indeed, we may
use Bing's house with the region $c$ marked as an o-spine of the
trivial knot. Whatever $r$, the only admissible coloring with color $0$
on the region $c$ has color $0$ on both of the other regions too. The weight of this
coloring is equal to $1$ whatever $r$ and $q$ are. The conclusion now follows by repeated
application of Theorem~\ref{ordinary:connected:thm}.
\end{proof}

The analogue of Corollary~\ref{eps:non-zero:color:cor} does not hold
for general~$r$. More precisely, we have the following:

\begin{prop}
There exist infinitely many pairs $r,q$ such that $\calI_2^{r,q}$ of the
unknot is not equal to 0.
\end{prop}

\begin{proof}
We use the abalone with the region $a$ marked to calculate the
invariants of the unknot. By definition we have
$$\calI_2^{r,q}(U_1)=
w_0w_2\left|\begin{array}{ccc}0 & 2 & 2 \\ 2 & 2 &
2\end{array}\right|+w_2^2\left|\begin{array}{ccc}2 & 2 & 2 \\ 2 &
2 & 2\end{array}\right|+w_2w_4\left|\begin{array}{ccc}4 & 2 & 2 \\
2 & 2 & 2\end{array}\right|\mbox{ for } r\geqslant 6.
$$
By direct calculation for all $r\geqslant 6$ we have
$$\calI_2^{r,q}(U_1)=1+[3]^2\frac{[5]-1}{[4]!}-\frac{1}{[4]},$$
and it is sufficient to show that for infinitely many $r$'s and some
choice of $q$ we have
$$[2][4]+[3][5]-[3]-[2]\neq 0.$$
Let us take $q=e^{(\pi\cdot\sqrt{-1})/r}$. By definition
we have
\begin{multline*}
[2][4]+[3][5]-[3]-[2]=\\
=\frac{q^8+q^{-8}+q^6+q^{-6}-(q^4+q^{-4})-
(q^2+q^{-2})-(q^3+q^{-3})+q+q^{-1}}{(q-q^{-1})^2}.
\end{multline*}
For the chosen $q$ this is equal to $0$ if and only if:
$$\cos\frac{8\pi}{r}+\cos\frac{6\pi}{r}-\cos\frac{4\pi}{r}-
\cos\frac{3\pi}{r}-\cos\frac{2\pi}{r}+\cos\frac{\pi}{r}=0.$$ If this
holds, $\cos\frac{\pi}{r}$ must be a solution of a specific equation
of degree $8$ with integer coefficients, namely the equation
$128t^8-224t^6+104t^4-4t^3-8t^2+4t=0$. Since this equation has at
most $8$ real roots, we get the conclusion.
\end{proof}

\paragraph{$\calI_0^{r,q}$ and the $TV$ invariant of the complement}
Here we show that, just as in the case of $\calI_0^{\varepsilon}$,
the invariants $\calI_0^{r,q}$ are determined by the link
complements. Namely, we obtain the following:

\begin{thm}\label{tv:link-compl:thm}
Let $L\subset M$ be a link in a closed $3$-manifold $M$, and let $N_M(L)$ be its
complement. Then
$$\calI_0^{r,q}(L)=TV(N_M(L)).$$
\end{thm}

The proof imitates that of Theorem~\ref{eps:link-compl:thm}, using Lemma~\ref{good-o-spine:lem}
and the analogue of Lemma~\ref{eps-0-puncture:lem}, which in turn is established
using \cite[Equation~(19)]{TV} in place of Lemma~\ref{eps:symbol-with-0:lem}.

\paragraph{Relations with the $\varepsilon$-invariant} It is
well-known that the $3$-manifold $\varepsilon$-invariant coincides
with the homologically trivial part of the Turaev-Viro invariant
of order $5$ for a suitable choice of $q$. It turns out that in the
case of knots in the sphere an even stronger statement holds.

\begin{thm}\label{eps:tv5:thm}
Let $K$ be a knot in $\matS^3$. Then
$$\calI_0^{\varepsilon}(K)=\calI_0^{5,q}(K)$$
for any choice of $\varepsilon$ and a suitable choice of $q$.
\end{thm}

\begin{lemma}
Let $K\subset\matS^3$ be a knot. Then for any $r$ and for any
projection $\Pi$ of $K$ any coloring in $\Adm_0(P_{\Pi}(K))$
involves integer colors only.
\end{lemma}

\begin{proof}
An admissible coloring on $P_{\Pi}(K)$ can be re-interpreted in the
following way. We assign a color to each region of the diagram
$\Pi$, to each overpass, to each underpass, to each crossing (this
is the color of the length-4 region arising from the tunnel going
over itself, see Fig.~\ref{crossing-point:fig}), and we choose two
more numbers, corresponding to the colors of the length-1 regions
arising with the ``puncturing'' move (see
Fig.~\ref{1-move-proj-spine:fig}). Notice that the underpass and the
overpass at the point where we add the transverse disc must
have the same color, and one of the regions adjacent to them must
have the same color as well. It follows that this color is an integer.
The region on the other side must also have an integer
color by admissibility.

It is easy to observe that any two regions separated by one edge of
the diagram must have colors that are either both integers or both
half-integers. Thus, since one region has integer color, the same
should be true for all other regions. Hence if an overpass and an
underpass intersect, their colors are either both integers or both
half-integers. Therefore the colors of all overpasses and of all
underpasses are either all integers or all half-integers; by the
above, they must be all integers.

It remains to check the conclusion for the regions corresponding to
crossing points and for the length-1 discs of
Fig.~\ref{1-move-proj-spine:fig}. But all these regions have
embedded closures, no two of them are adjacent to the same edge and
all the other regions have already been shown to have integer
colors. Then by admissibility all these regions must indeed have
integer colors.
\end{proof}\vspace{3mm}

\noindent\emph{Proof of Theorem~\ref{eps:tv5:thm}.} Notice that any
admissible coloring in the sense of the initial datum
$\mathcal{D}_{\varepsilon}$ is also an admissible coloring
in the sense of any Turaev-Viro initial datum of order 5. If
we choose $q$ to be such that $\varepsilon=q+q^{-1}$, which we can always do,
this correspondence preserves the weights
of the colors and the values of quantum $6j$-symbols, whence the
conclusion. \finedimo\vspace{5mm}

We remark that calculation of $\calI_0^{\varepsilon}$ is much easier
than the calculation of $\calI_0^{5,q}$, hence for knots it is more
effective to study $\calI_0^{\varepsilon}$ rather than $\calI_0^{5,q}$.
On the other hand, there exist links for which
Theorem~\ref{eps:tv5:thm} does not hold.

\section{Relations with polynomial
invariants}\label{other-invts:section}

In this section we consider the relation between the colored Turaev-Viro
invariants and some of the well-known polynomial invariants of
links, namely, the HOMFLY polynomial, the Kauffman polynomial, and
the Alexander polynomial.

\begin{thm}\label{eps:distinct:skein:thm}
There exists a pair of links in $\matS^3$ such that their HOMFLY
polynomials are the same and their Kauffman polynomials are the
same but their $\varepsilon$-invariants $\calI_0^{\varepsilon}$
are distinct.
\end{thm}

\begin{proof}
Recall that both the Kauffman polynomial $F_L(a,z)$ and the HOMFLY
polynomial $P_L(a,z)$ are multiplicative under connected sum of
links. In particular, let $L=K_1\sqcup K_2$ be a two-component link,
and let $K$ be a knot. Denote by $L^{(i)}$ the link
obtained by taking the connected sum of $L$ with $K$ along $K_i$.
Then:
$$\begin{array}{c}
F_{L^{(1)}}(a,z)=F_L(a,z)\cdot F_K(a,z)=F_{L^{(2)}}(a,z),\\
P_{L^{(1)}}(a,z)=P_L(a,z)\cdot P_K(a,z)=P_{L^{(2)}}(a,z).
\end{array}$$

Let us give an example of $L$ and $K$ such that
$\calI_0^{\varepsilon}(L^{(1)})\neq\calI_0^{\varepsilon}(L^{(2)})$.
Take $K$ to be a trefoil and take $L$ to be the link shown in
Fig.~\ref{link-tref-triv:fig}.
\begin{figure}
    \begin{center}
    \includegraphics[scale=0.5]{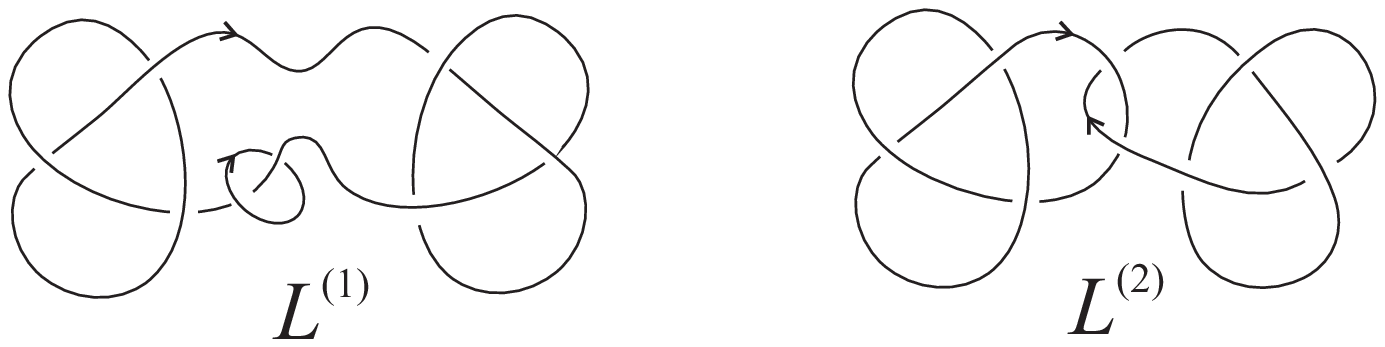}
    \mycap{The links $L^{(1)}$ and $L^{(2)}$}
    \label{two-links-from-tref-triv:fig}
    \end{center}
    \end{figure}
Denote by $K_1$ the trefoil component of $L$ and by $K_2$ the
trivial component of $L$. The links $L^{(1)}$ and $L^{(2)}$ are
shown in Fig.~\ref{two-links-from-tref-triv:fig}. Then by
Theorem~\ref{eps:2-connected:thm} and Corollary~\ref{eps:unlink:cor} we have:
\begin{multline*}
\calI_0^{\varepsilon}(L^{(1)})=
\calI_0^{\varepsilon}(L)\cdot\calI_0^{\varepsilon}(3_1)\\
+\varepsilon^{-2}\left((\varepsilon+2)\calI_0^{\varepsilon}(U_1)-\calI_0^{\varepsilon}(L)
\right)\cdot\left((\varepsilon+2)\calI_0^{\varepsilon}(\emptyset)-\calI_0^{\varepsilon}(3_1)
\right)=18-7\varepsilon.
\end{multline*}
At the same time we have
\begin{multline*}
\calI_0^{\varepsilon}(L^{(2)})=
\calI_0^{\varepsilon}(L)\cdot\calI_0^{\varepsilon}(3_1)\\
+\varepsilon^{-2}\left((\varepsilon+2)\calI_0^{\varepsilon}(3_1)-\calI_0^{\varepsilon}(L)
\right)\cdot\left((\varepsilon+2)\calI_0^{\varepsilon}(\emptyset)-\calI_0^{\varepsilon}(3_1)
\right)=15-6\varepsilon,
\end{multline*}
whence the conclusion.
\end{proof}

\begin{rem}
\emph{The links constructed in the proof of the above theorem are
reducible, but there also exists pairs of prime knots with the
same HOMFLY polynomial but distinct Turaev-Viro invariants of
order $7$. One such example is given by the knots $8_8$ and
$10_{129}$ in the tables of Tait~\cite{Tait} and
Little~\cite{Little} (see Fig.~\ref{irr-same-homfly:fig}), which
have the same HOMFLY polynomial~\cite{Jones}.
\begin{figure}
    \begin{center}
    \includegraphics[scale=0.5]{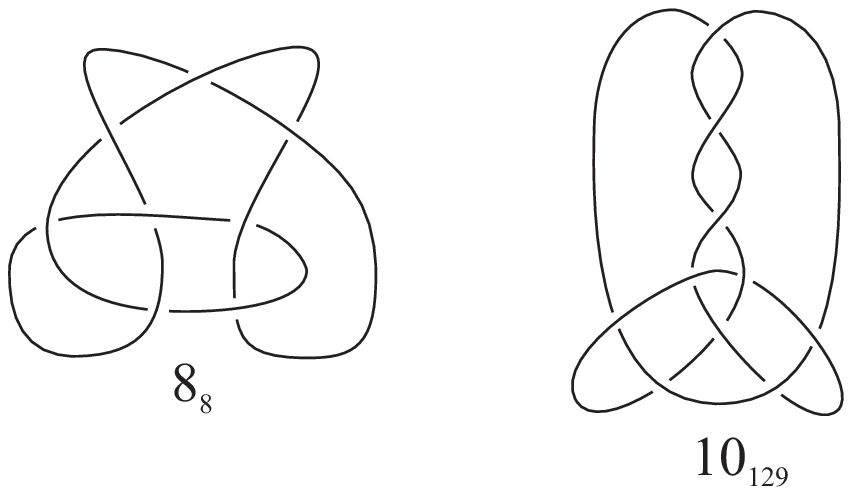}
    \mycap{Prime knots with the same HOMFLY polynomial and distinct
    Turaev-Viro invariants}
    \label{irr-same-homfly:fig}
    \end{center}
    \end{figure}
We have calculated their invariants $\calI_0^{7,q}$ and $\calI_0^{8,q}$ using
Theorem~\ref{eps:link-compl:thm} and the program ``Three-manifold
Recognizer'' \cite{MathRec}, which, among other things, calculates
Turaev-Viro invariants of link complements in $\matS^3$. The results
were:
$$
\left.\begin{array}{lcl} \calI_0^{7,q}(8_8) & = &
71(q^2+q^{-2})-181(q+q^{-1})+240,\\
\calI_0^{7,q}(10_{129}) & = & 89(q^2+q^{-2})-227(q+q^{-1})+301,\\
 & & \\
\calI_0^{8,q}(8_8) & = &
79-54(q^2+q^{-2}),\\
\calI_0^{8,q}(10_{129}) & = & 31-20(q^2+q^{-2}).
\end{array}\right.
$$ }
\end{rem}

\begin{prop}\label{eps:distinct:Alex:thm}
There exists a pair of knots with the same Alexander polynomial
but distinct $\varepsilon$-invariants.
\end{prop}

\begin{proof}
Let $W(K)$ denote the Whitehead double of a given knot
$K\subset\matS^3$. Then it is known
that the Alexander polynomial of
$W(K)$ is equal to $1$. On the other hand, a direct (but computer-aided) calculation
allowed us to obtain that
$$\calI_0^{\varepsilon}(W(3_1))=5-2\varepsilon,\;\;
\calI_0^{\varepsilon}(W(4_1))=4-\varepsilon. $$
\end{proof}

\section{Colored Turaev-Viro invariants\\ not determined by
link complements}\label{invt-not-complement:sect}

Theorems~\ref{eps:link-compl:thm} and~\ref{tv:link-compl:thm} show that some
colored Turaev-Viro invariants of a link are determined by the link
complement. We now show that not all of them are.

\begin{prop}
There exists a pair of links in $\matS^3$ with homeomorphic
complements but distinct $\calI_1^{\varepsilon}$ invariants.
\end{prop}

\begin{proof}
Consider the links obtained by closing
the braids $\sigma_2\sigma_1^2\sigma_2$ and
$\sigma_1^2\sigma_2\sigma_1^2\sigma_2$, see
Fig.~\ref{same-compl-links:fig}.
\begin{figure}
    \begin{center}
    \includegraphics[scale=0.3]{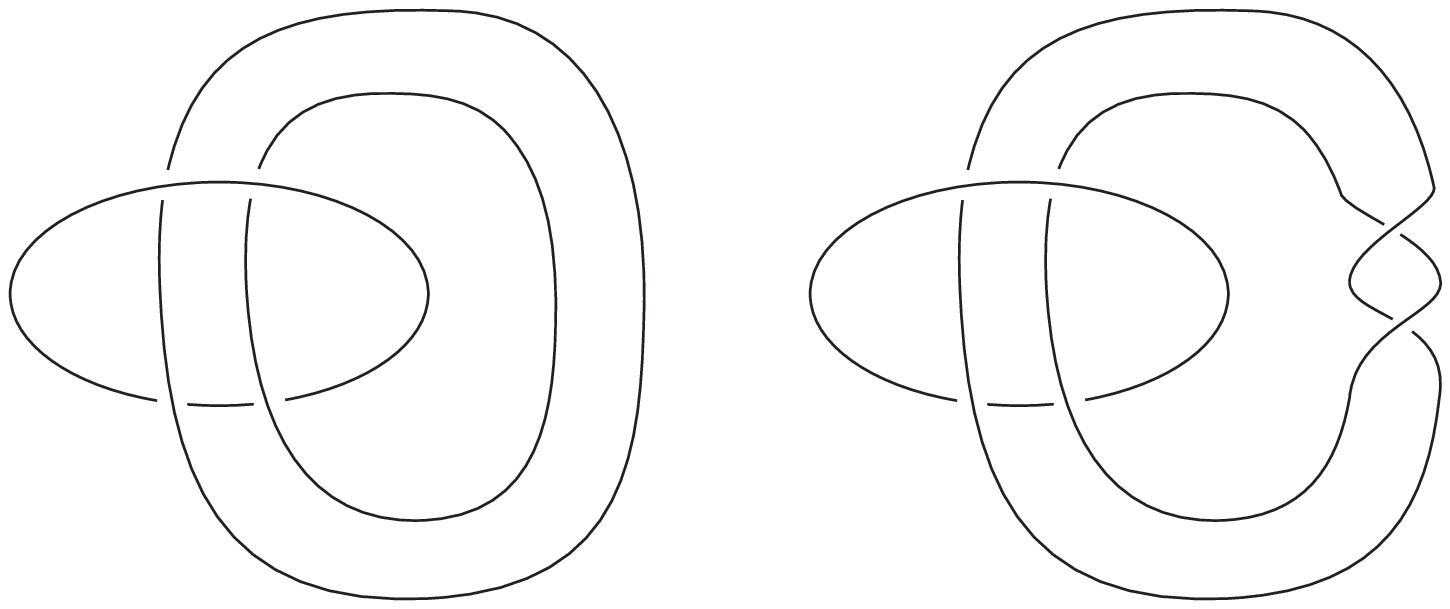}
    \mycap{Two links with the same complement}
    \label{same-compl-links:fig}
    \end{center}
    \end{figure}
Since they are obtained from each other by a full twist around a
disc spanning the common unknotted component, they have the same complements.
To compute their $\calI_1^{\varepsilon}$, we constructed o-spines from the
projections, obtaining the polyhedra of
Fig.~\ref{same-compl-ospines:fig}.
\begin{figure}
    \begin{center}
    \includegraphics[scale=0.4]{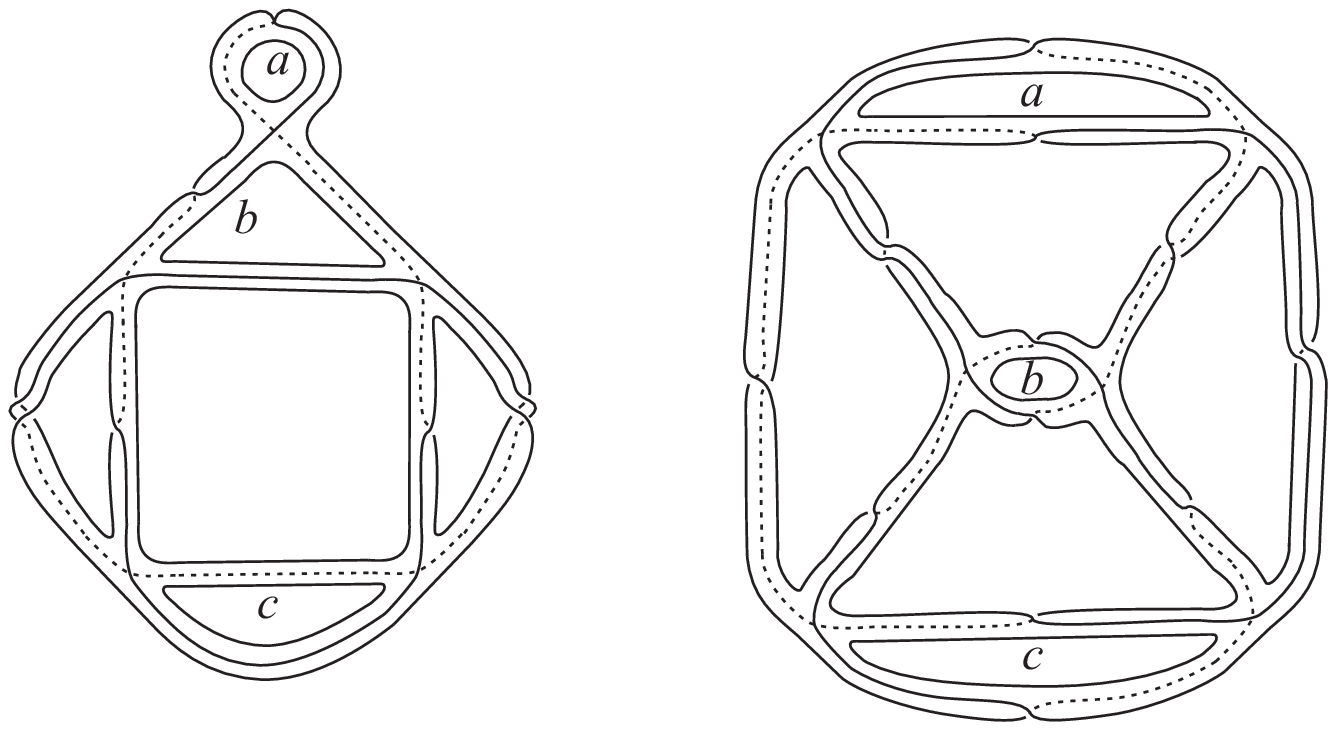}
    \mycap{O-spines of the links of Fig.~\ref{same-compl-links:fig}.
    The letters $a$, $b$, $c$ denote the regions dual to the link components}
    \label{same-compl-ospines:fig}
    \end{center}
    \end{figure}
A direct computation then shows that
$$\calI_1^{\varepsilon}(\widehat{\sigma_2\sigma_1^2\sigma_2})=-(\varepsilon+2),$$
$$\calI_1^{\varepsilon}(\widehat{\sigma_1^2\sigma_2\sigma_1^2\sigma_2})=-(4\varepsilon+3).$$
\end{proof}

\noindent
Dipartimento di Matematica Applicata\\
Universit\`a di Pisa\\
Via Filippo Buonarroti, 1C\\
I-56127, PISA\\
Italy\\
\ \\
{\tt pervova@guest.dma.unipi.it}\\
{\tt petronio@dm.unipi.it}

\end{document}